\documentclass[a4paper]{scrartcl}
\usepackage{mathtools}
\usepackage{amssymb}
\usepackage{graphicx}
\usepackage{subfigure}
\usepackage{url}
\usepackage{algorithm}
\usepackage[noend]{algpseudocode}
\usepackage{tikz}
\usetikzlibrary{patterns}

\newcommand{\cU}{\mathcal{U}}
\newcommand{\X}{\mathcal{X}}
\newcommand{\R}{\mathbb{R}}
\newcommand{\cR}{\mathcal{R}}

\definecolor{vgreen}{HTML}{188F00}
\definecolor{vblue}{HTML}{3500AE}

\makeatletter
\def\BState{\State\hskip-\ALG@thistlm}
\makeatother

\usepackage{authblk}

\begin{document}

% Author macros::begin %%%%%%%%%%%%%%%%%%%%%%%%%%%%%%%%%%%%%%%%%%%%%%%%
% \title{An Experimental Comparison of Uncertainty Sets for Robust Shortest Path Problems}
\title{Algorithms and Uncertainty Sets for Data-Driven Robust Shortest Path Problems}
% \titlerunning{A Sample LIPIcs Article} %optional, in case that the title is too long; the running title should fit into the top page column

\author[1]{Andr\'{e} Chassein\thanks{Email: chassein@mathematik.uni-kl.de}}

\author[2]{Trivikram Dokka\thanks{Email: t.dokka@lancaster.ac.uk}}

\author[2]{Marc Goerigk\thanks{Corresponding author. Email: m.goerigk@lancaster.ac.uk}}

\date{}

\affil[1]{Department of Mathematics, University of Kaiserslautern, Kaiserslautern, Germany}
\affil[2]{Department of Management Science, Lancaster University, Lancaster, United Kingdom}

\maketitle

\begin{abstract}
We consider robust shortest path problems, where the aim is to find a path that optimizes the worst-case performance over an uncertainty set containing all relevant scenarios for arc costs. The usual approach for such problems is to assume this uncertainty set given by an expert who can advise on the shape and size of the set.

Following the idea of data-driven robust optimization, we instead construct a range of uncertainty sets from the current literature based on real-world traffic measurements provided by the City of Chicago. We then compare the performance of the resulting robust paths within and outside the sample, which allows us to draw conclusions what the most suited uncertainty set is.

Based on our experiments, we then focus on ellipsoidal uncertainty sets, and develop a new solution algorithm that significantly outperforms a state-of-the-art solver.
\end{abstract}

\textbf{Keywords:} robustness and sensitivity analysis; robust shortest paths; uncertainty sets; data-driven robust optimization

\section{Introduction}

For classic shortest path problems in street networks, considerable speed-ups over a standard Dijkstra's algorithm have been achieved thanks to algorithm engineering techniques \cite{bast2016route}, which makes real-time information even in large networks possible. Most types of robust shortest path problems, on the other hand, are NP-hard (see \cite{yu1998robust}), and real-time information has not been an option.

To formulate a robust problem, it is necessary to have a description of all possible and relevant scenarios that the solution should prepare against, the so-called uncertainty set. We refer to the surveys  \cite{aissi2009min,goerigk2016algorithm,kasperski2016robust} for a general overview on the topic. The current literature on robust shortest paths usually assumes this set to be given, by some mixture of data-preprocessing and expert knowledge that is not part of the study. This means that different types of sets have been studied (compare, e.g., \cite{montemanni2004exact,busing2012recoverable}), but it has been impossible to address the question which would be the ''right'' choice.

A recent paradigm shift is data-driven robust optimization (see \cite{bertsimas2017data}), where building the uncertainty set from raw observations is part of the robust optimization problem. This paper is the first to follow such a perspective for shortest path problems. Based on real-world observations from the City of Chicago, we build a range of uncertainty sets, calculate the corresponding robust solutions, and perform an in-depth analysis of their performance. This allows us to give an indication which set is actually suitable for our application, and which are not. 

In the second part of this paper, we then focus on the case of ellipsoidal uncertainty, and provide a branch-and-bound algorithm that is able to solve instances considerably faster than an off-the-shelf solver.

Parts of this paper were previously published as a conference paper in \cite{atmos}. In comparison, we provide a completely new set of experimental results based on an observation period of 46 days (instead of one single day), which leads to a more detailed insight into the performance of different uncertainty sets. Furthermore, we provide a new analysis for axis-parallel ellipsoidal uncertainty sets, including an efficient branch-and-bound algorithm that is able to outperform Cplex by several orders of magnitude, pushing robust shortest paths towards applicability in real-time navigation systems.

The remainder of the paper is structured as follows. In Section~\ref{sec:unc} we briefly introduce the robust shortest path problem along with six uncertainty sets used in this study. The experimental setup and results on real-world data are then presented in Section~\ref{sec:exp}. Section~\ref{sec:ellipsoids} presents our algorithm for ellipsoidal uncertainty sets, and includes additional computational results. The paper is concluded in Section~\ref{sec:conclusion}.

\section{Uncertainty Sets for the Shortest Path Problem}
\label{sec:unc}

In this section, we briefly introduce the robust shortest path problem and different approaches to model uncertainty sets that are used in the current literature. In the classic shortest path problem, we are given a directed graph $G=(V,A)$ where $V$ denotes the set of nodes, and $A$ denotes the set of arcs. For every arc $e\in A$, we know its traversal time $c_e \ge 0$. For a start node~$s$ and a target node~$t$, the aim is to find a path minimizing the total travel time given as the sum of times over all arcs that are part of the path.

More formally, we denote this problem as 
\[ \min \left\{ \pmb{c}^t\pmb{x} : \pmb{x}\in\X \right\} \]
$\X\subseteq\{0,1\}^n$ denotes the set of $s$-$t$-paths, and $n = |A|$ is the number of variables.

\begin{figure}[htbp]
\centering
\subfigure[Convex hull with \textcolor{vblue}{$\lambda=1$} and \textcolor{vgreen}{$\lambda=0.5$}.]{\label{graph1}\includegraphics[width=0.45\textwidth]{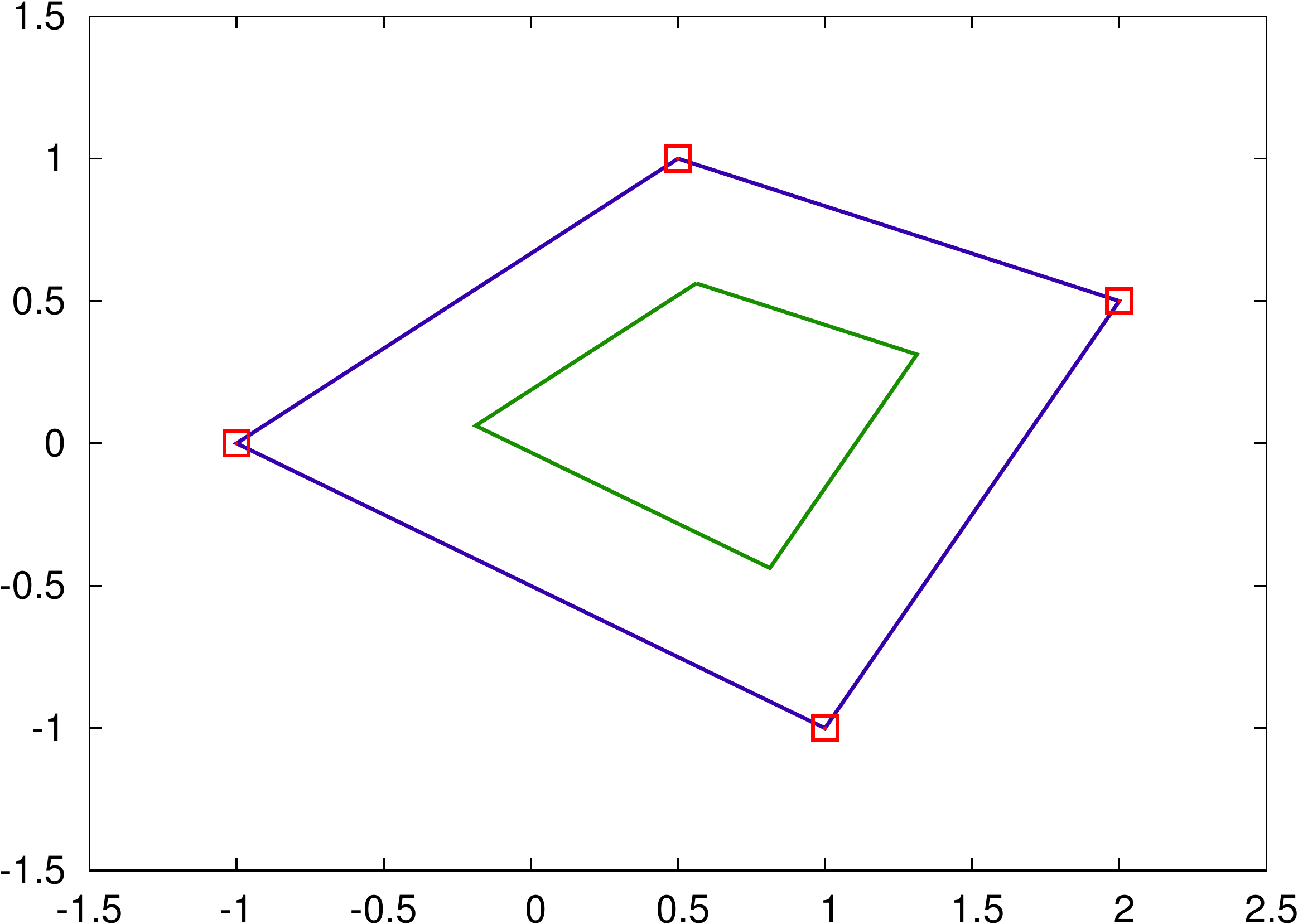}}  
    \hfill
\subfigure[Intervals with \textcolor{vblue}{$\lambda=1$} and \textcolor{vgreen}{$\lambda=0.5$}.]{\label{graph2}\includegraphics[width=0.45\textwidth]{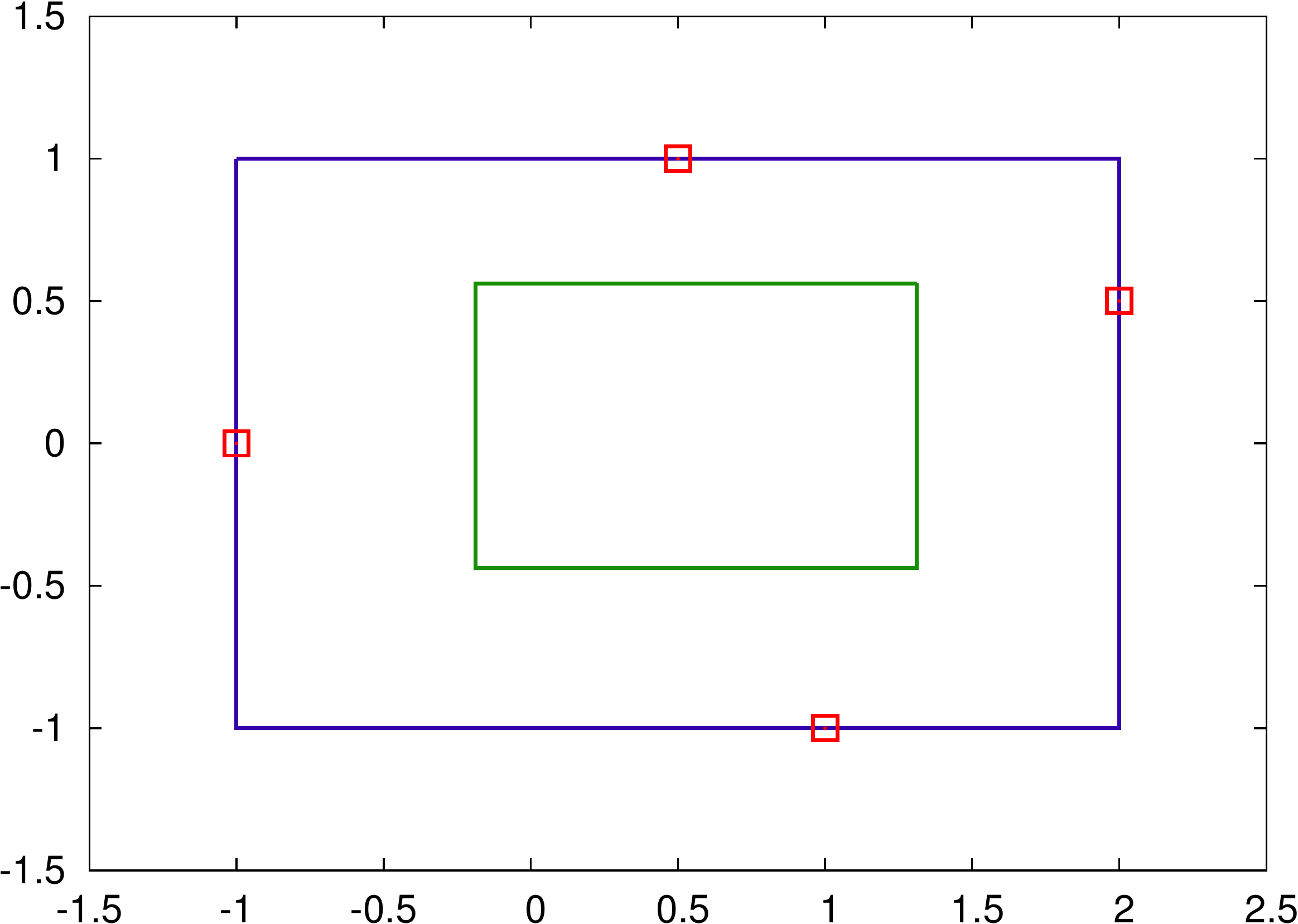}}
\subfigure[Ellipsoid with \textcolor{vblue}{$\lambda=3$} and \textcolor{vgreen}{$\lambda=1$}.]{\label{graph3}\includegraphics[width=0.45\textwidth]{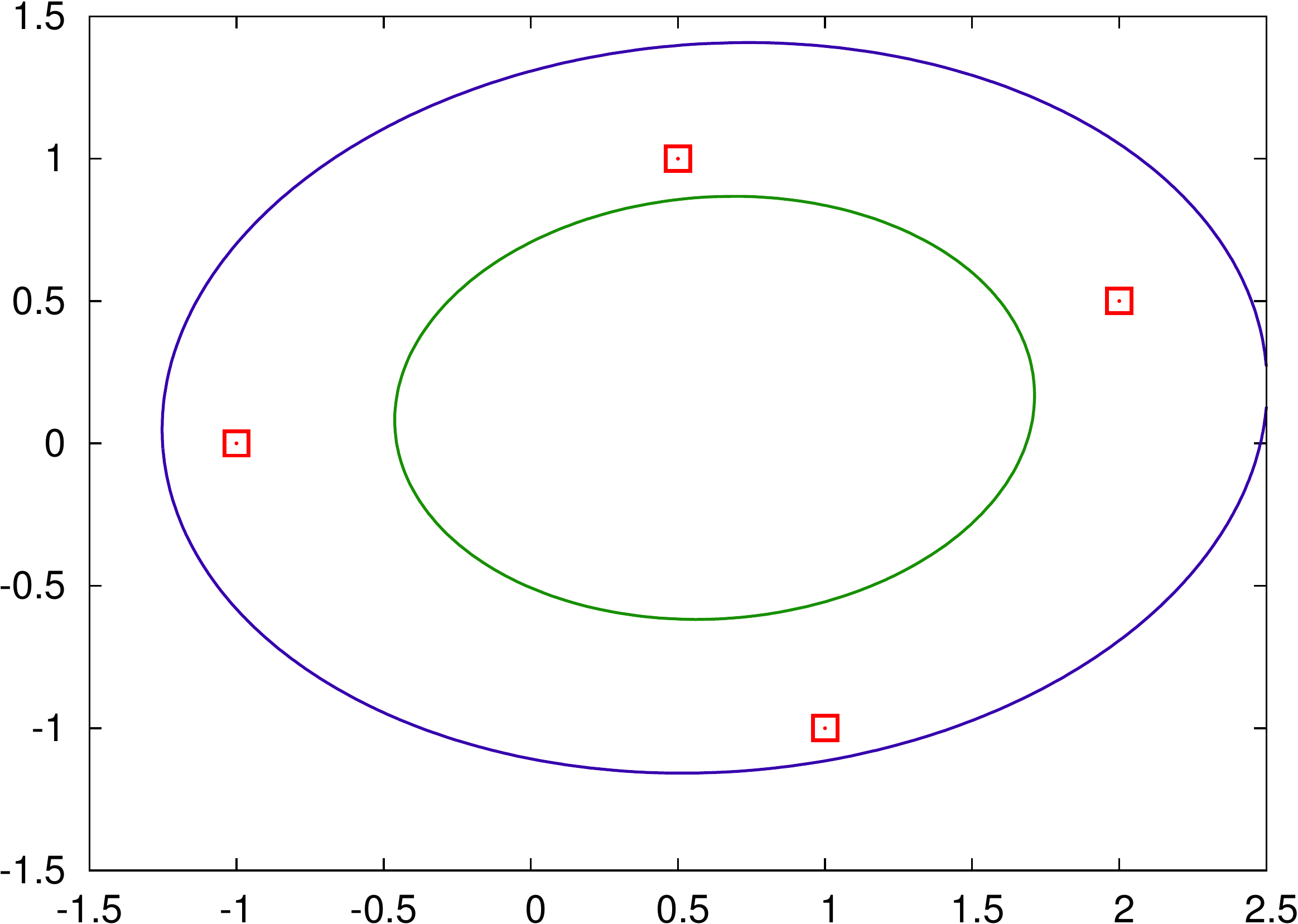}}      
    \hfill
\subfigure[Budgeted uncertainty with \textcolor{vblue}{$\Gamma=1.5$} and \textcolor{vgreen}{$\Gamma=1$}.]{\label{graph4}\includegraphics[width=0.45\textwidth]{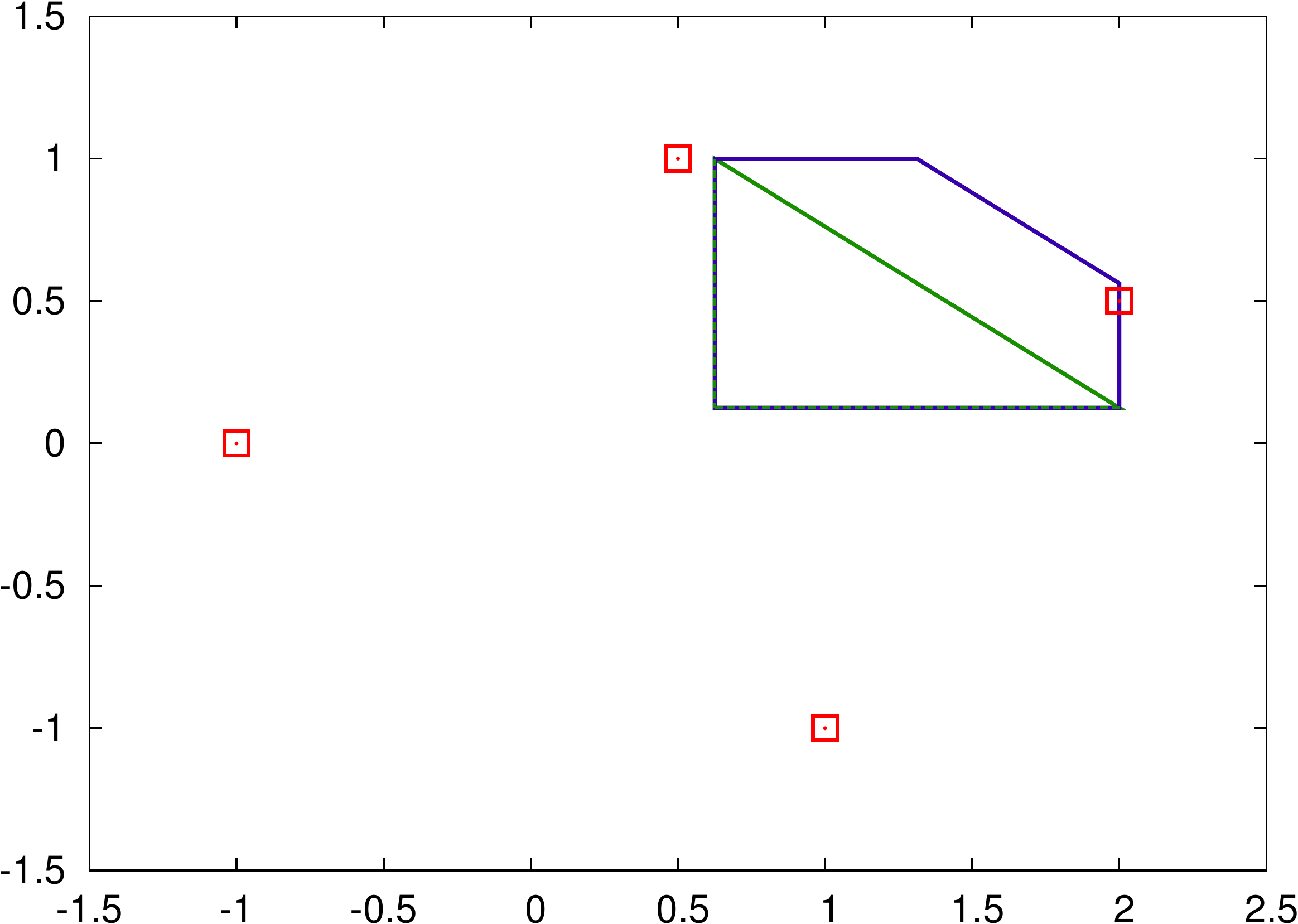}}
\subfigure[Permutohull uncertainty for \textcolor{vblue}{$CVaR_{2/N}$} and \textcolor{vgreen}{$CVaR_{3/N}$} (i.e., \textcolor{vblue}{$\pmb{q} = (\frac{1}{2},\frac{1}{2},0,0)$} and \textcolor{vgreen}{$\pmb{q}=(\frac{1}{3},\frac{1}{3},\frac{1}{3},0)$}.]{\label{graph5}\includegraphics[width=0.45\textwidth]{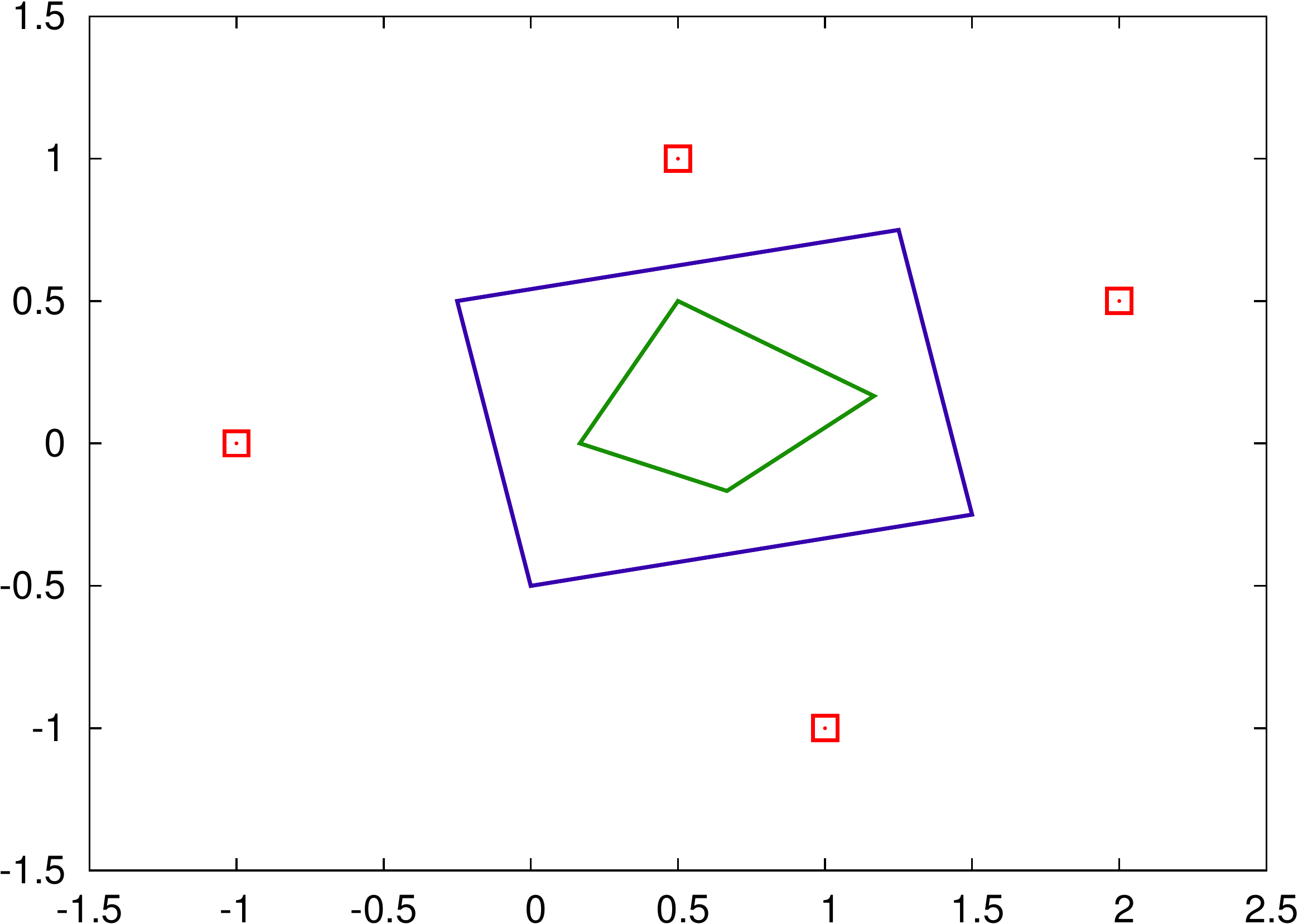}}
    \hfill
\subfigure[Symmetric permutohull uncertainty for \textcolor{vblue}{$\pmb{q}=(\frac{1}{2},\frac{1}{2},0,0)$} and \textcolor{vgreen}{$\pmb{q}=(\frac{1}{2},\frac{1}{4},\frac{1}{4},0)$}.]{\label{graph6}\includegraphics[width=0.45\textwidth]{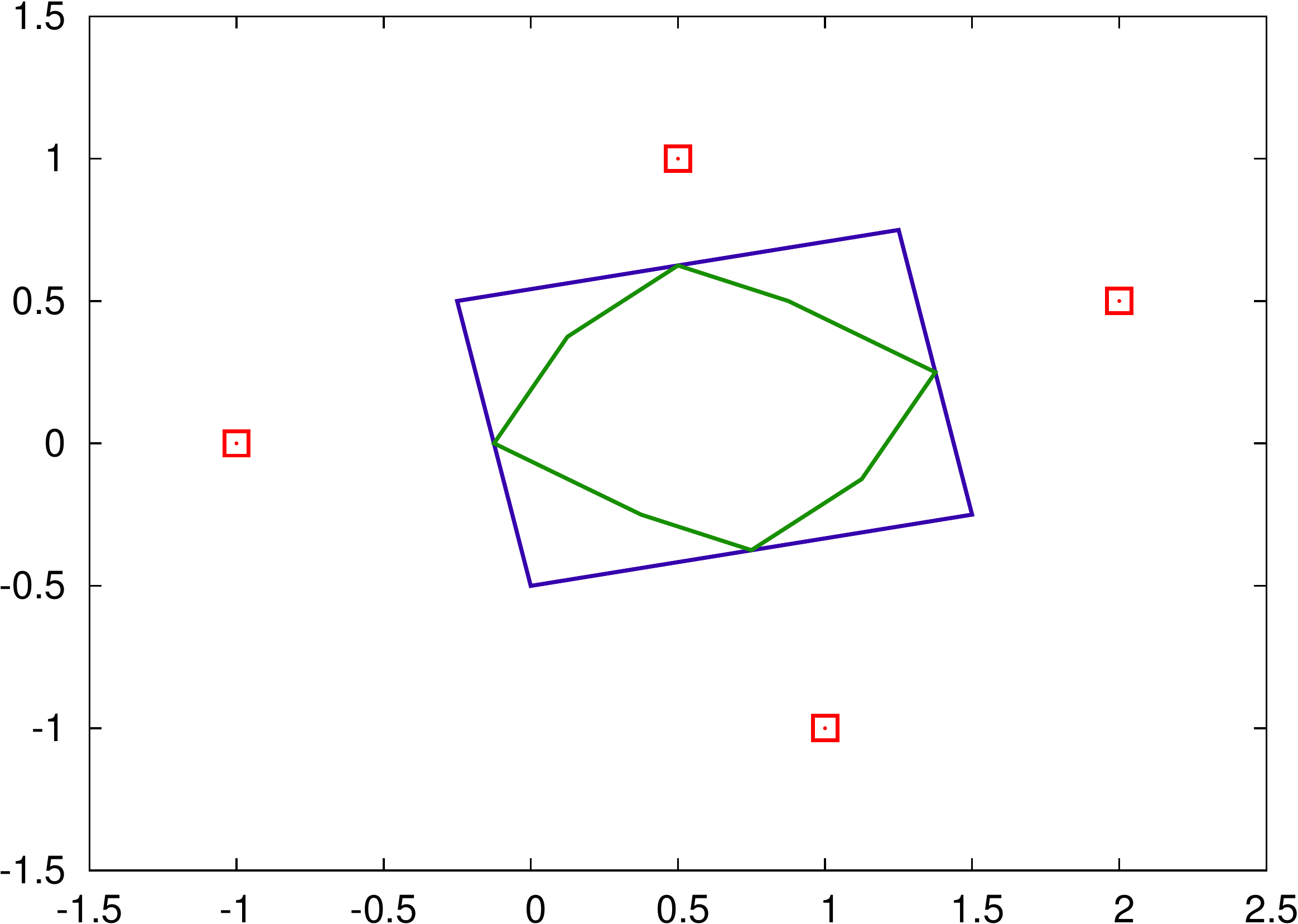}}
\caption{Example uncertainty sets.}\label{fig1}
\end{figure}

In our setting, we assume that travel times $\pmb{c}$ are not known exactly. Instead, we are provided with a set $\cR$ of travel time observations, where $\cR = \{\pmb{c}^1,\ldots,\pmb{c}^N\}$ with $\pmb{c}^i \in \mathbb{R}^n$. We also refer to $\cR$ as the available raw data. Based on this data, an uncertainty set $\cU$ is generated which is then used within the robust shortest path problem
\[ \min \left\{ \max_{\pmb{c}\in\cU} \pmb{c}^t\pmb{x} : \pmb{x}\in\X \right\} \]
that is, we search for a path that minimizes the worst-case costs over all costs in $\cU$. We briefly sketch approaches to generate $\cU$ in the following. Each is equipped with a scaling parameter to control its size (see also \cite{chassein2017variable} on the problem of choosing the size of an uncertainty set with a given shape). A visual example using four data points in two dimensions is provided for each apprach in Figure~\ref{fig1}. We use the notation $[N]=\{1,\ldots,N\}$ and denote by $\hat{\pmb{c}}$ the average of $\{\pmb{c}^1,\ldots,\pmb{c}^N\}$, i.e., $\hat{\pmb{c}} = \frac{1}{N}\sum_{i\in[N]} \pmb{c}^i$. For more details on the resulting models, we refer to the conference version of this paper \cite{atmos}.

\begin{itemize}
\item Convex hull uncertainty (see, e.g., \cite{kasperski2016robust,yu1998robust}): We set 
\[ \cU^{CH} = \cR.\]
Note that this is equivalent to using the convex hull of raw data. To scale this set, we substitute each point $\pmb{c}^i$ with $\hat{\pmb{c}} + \lambda(\pmb{c}^i-\hat{\pmb{c}})$ for a given $\lambda\ge 0$, and take the convex hull of the scaled data points.

\item Interval uncertainty (see, e.g., \cite{chassein2015new}): We set 
\[ \cU^I = \bigtimes_{i\in[n]} \left[ \hat{c}_i + \lambda (\underline{c}_i - \hat{c}_i), \hat{c}_i + \lambda (\overline{c}_i - \hat{c}_i) \right]\]
for some $\lambda\ge0$.

\item Ellipsoidal uncertainty (see, e.g., \cite{ben1998robust,ben1999robust}): We set 
\[ \cU^E = \{ \pmb{c} : (\pmb{c} - \pmb{\mu})^t\pmb{\Sigma}^{-1}(\pmb{c} - \pmb{\mu}) \le \lambda \}\]
with $\pmb{\mu} = \hat{\pmb{c}} = \frac{1}{N}(\pmb{c}^1 + \ldots + \pmb{c}^N)$ and $\pmb{\Sigma} = \frac{1}{N} \sum_{i\in[N]} (\pmb{c}^i - \pmb{\mu})(\pmb{c}^i - \pmb{\mu})^t$ derived from a maximum-likelihood fit of a normal distribution.

\item Budgeted uncertainty (see, e.g., \cite{bertsimas2003robust,bertsimas2004price,goerigk2016algorithm}): We set \[\cU^{B} = \{ \pmb{c} : c_i = \hat{c}_i + (\overline{c}_i - \hat{c}_i)\delta_i \text{ for all } i\in[n],\ \pmb{0} \le \pmb{\delta}\le \pmb{1},\ \sum_{i\in[n]} \delta_i \le \Gamma\}\]
where the parameter $\Gamma$ controls the size of $\cU^{B}$.

\item Permutohull uncertainty (see \cite{bertsimas2009constructing}): We set
\[ \cU^{PH} = conv\left( \left\{ \sum_{i\in[N]} q_{\sigma(i)} \pmb{c}^i : \sigma \in S_N\right\}\right) \]
where $S_N$ denotes the set of permutations on $[N]$, and $\pmb{q}$ is a column of the matrix
\[ Q_N := \begin{pmatrix}
1 & \dots & \frac{1}{N-2} & \frac{1}{N-1} & \frac{1}{N} \\
\vdots & \vdots & \vdots & \vdots & \vdots \\
0 & 0 & \frac{1}{N-2} & \frac{1}{N-1} & \frac{1}{N} \\ 
0 & \dots & 0 & \frac{1}{N-1} & \frac{1}{N} \\
0 & \dots & 0 & 0 & \frac{1}{N}
\end{pmatrix} \in \R^{N\times N}\]
Scaling is included via the choice of the column, where using the $j$th column of $Q_N$ corresponds to using the conditional value at risk $CVaR$ criterion with respect to risk level $j/N$.

\item Symmetric permutohull uncertainty (see \cite{bertsimas2009constructing}): As in the above case, but we generate $\cU^{SPH}$ using columns of the matrix $\tilde{Q}\in\R^{N \times (\lfloor N/2 \rfloor +1)}$ defined by
\[ \tilde{Q} := \frac{1}{N}\begin{pmatrix}
1 & 2 & 2 & \dots & 2 \\
1 & 1 & 2 & \dots & 2 \\
1 & 1 & 1 & \dots & 2 \\
\vdots & \vdots & \vdots  & \vdots  & \vdots  \\
1 & 1 & 1 & \dots & 0 \\
1 & 1 & 0 & \dots & 0 \\
1 & 0 & 0 & \dots & 0 \\
\end{pmatrix}, \]
instead, i.e., in the first column, all entries are $1/N$; in the second column, the first entry is $2/N$ and the last entry is $0$, etc. The resulting sets are symmetric with respect to $\hat{\pmb{c}}$.

\end{itemize}
 
In total we use six methods to generate uncertainty set $\cU$ based on the raw data $\cR$. The resulting optimization model and its complexity are summarized in Table~\ref{tab1}. Here, "(M)IP" stands for (mixed-)integer linear program, "LP" for linear program, and "MISOCP" for mixed-integer second order cone program.
\begin{table}[htb]
\centering\begin{tabular}{r|c|c|c|c|c|c}
 & $\cU^{CH}$ & $\cU^{I}$ & $\cU^E$ & $\cU^{B}$ & $\cU^{PH}$ & $\cU^{SPH}$ \\
 \hline
Complexity & NPH & P & NPH & P & NPH & NPH \\
Model & IP & LP & MISOCP & MIP & MIP & MIP \\
Add. Const. & $N$ & 0 & 1 & $n+1$ & $N^2$ & $N^2$ \\
Add. Var. & 1 & 0 & 1 & $n$ & $2n$ & $2n$ \\
\end{tabular}
\caption{Uncertainty sets in this study.}\label{tab1}
\end{table}
While the robust model with budgeted uncertainty sets can be solved in polynomial time using combinatorial algorithms, we still used the MIP formulation for our experiments, as it was sufficiently fast.

\section{Real-World Experiments}
\label{sec:exp}
 
\subsection{Data Collection and Cleaning}

We used data provided by the City of Chicago\footnote{https://data.cityofchicago.org}, which provides a live traffic data interface. The data set consists of traffic updates for every 15-minute interval over a time horizon of 46 days spanning Tuesday morning, March 28, 2017 to Friday evening, May 12, 2017.

\begin{figure}[htbp]
\centering
\subfigure[Raw segments with zoom-in for the city center, longitude versus latitude. In red are segments without data.]{\label{chicago1}\includegraphics[width=.7\textwidth]{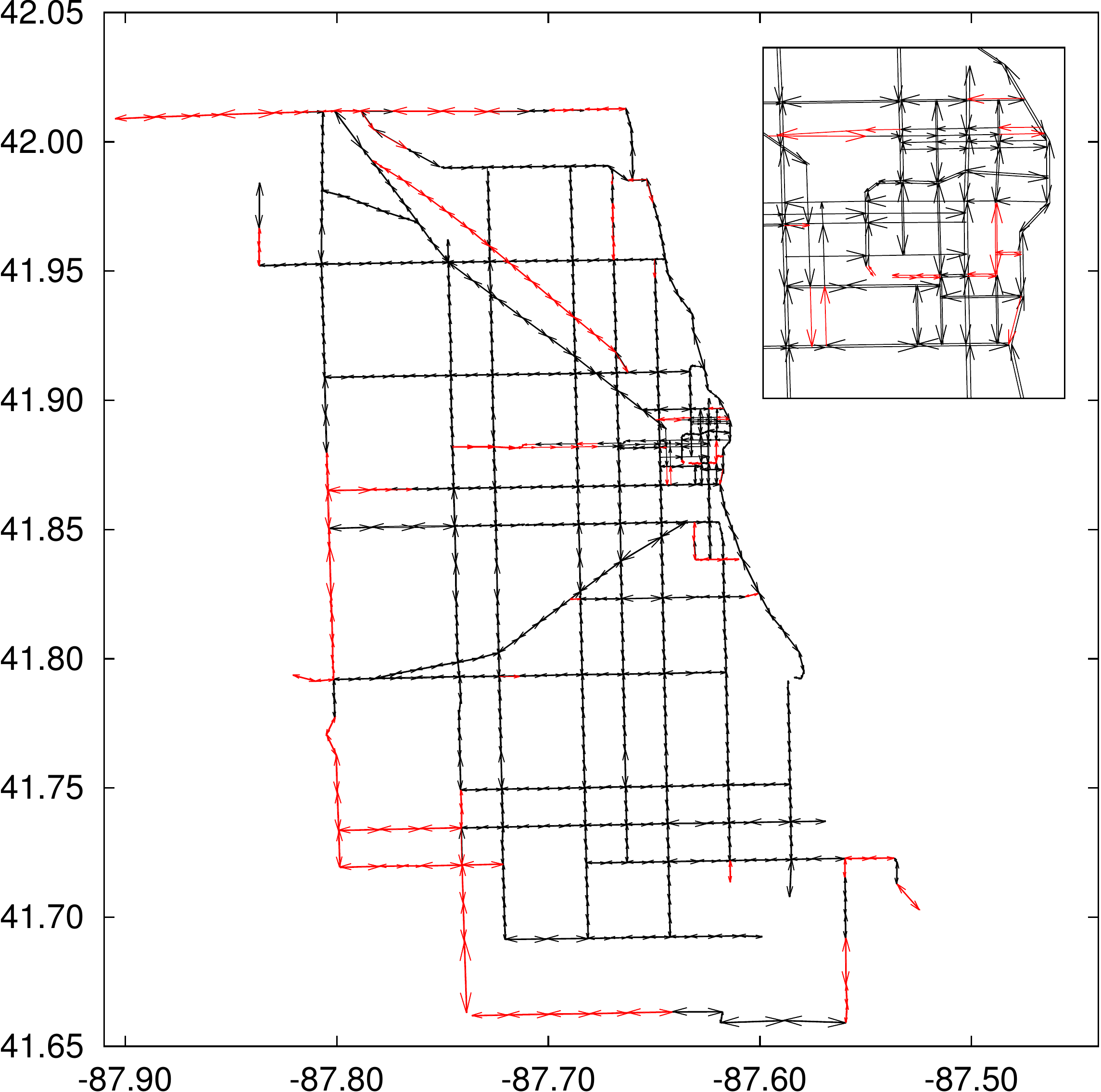}}  
    \hfill
\subfigure[Resulting graph model with zoom-in for the city center, longitude versus latitude.]{\label{chicago2}\includegraphics[width=.7\textwidth]{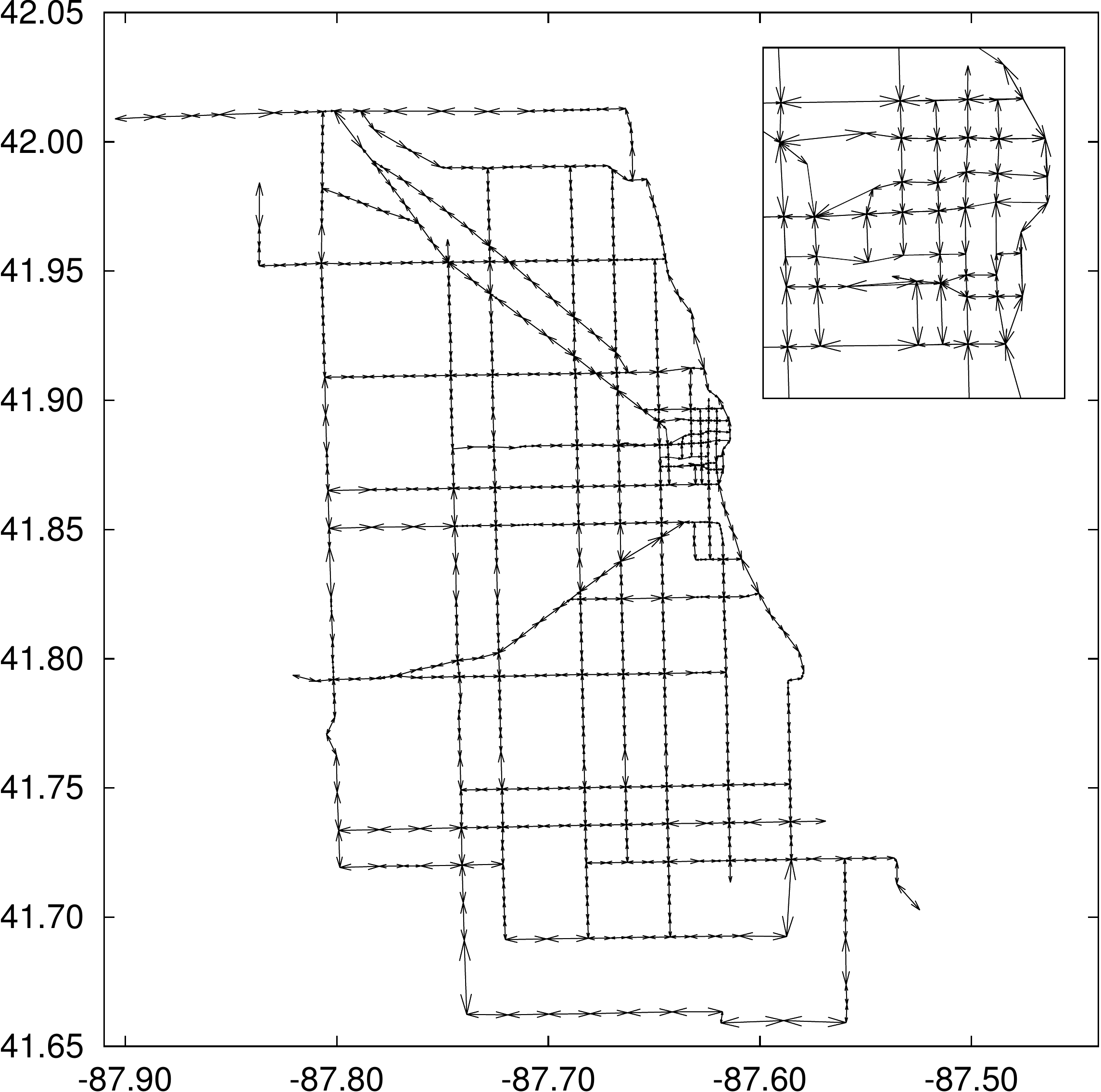}}     
\caption{Chicago instance.}
\end{figure}

Out of all $46\cdot 96 = 4416$ potential observations, only 4363 had usable data or were recorded due to server downtimes. Every data point contains the traffic speed for a subset of a total of 1,257 segments. For each segment the geographical position is available, see the resulting plot in Figure~\ref{chicago1} with a zoom-in for the city center. 
The complete travel speed data set contains a total of 3,891,396 records. There were 1,045 segments where the data was recorded at least once in the 4363 data points. 
For nearly 55\% of the segments, at least 1340 data points were recorded, and more than 90\% of them have at least 450 data points. For almost all segments at least 400 data points were recorded.
We used linear interpolation to fill the missing records keeping in mind that data was collected over time.
For segments that did not have any data, we set the travel speed to 20 miles per hour (which is slightly slower than the average speed in the network). Any speed record below 3 miles per hour was set to 3 miles per hour to ensure resonable travel times.
Segment lengths were given through longitude and latitude coordinates, and approximated using the Euclidean distance.

% The data after removing missing records and filling missing values can be found at \url{www.lancaster.ac.uk/~goerigk/robust-sp-data.zip}. 

\begin{figure}[htbp]
\centering
\includegraphics[width=\textwidth]{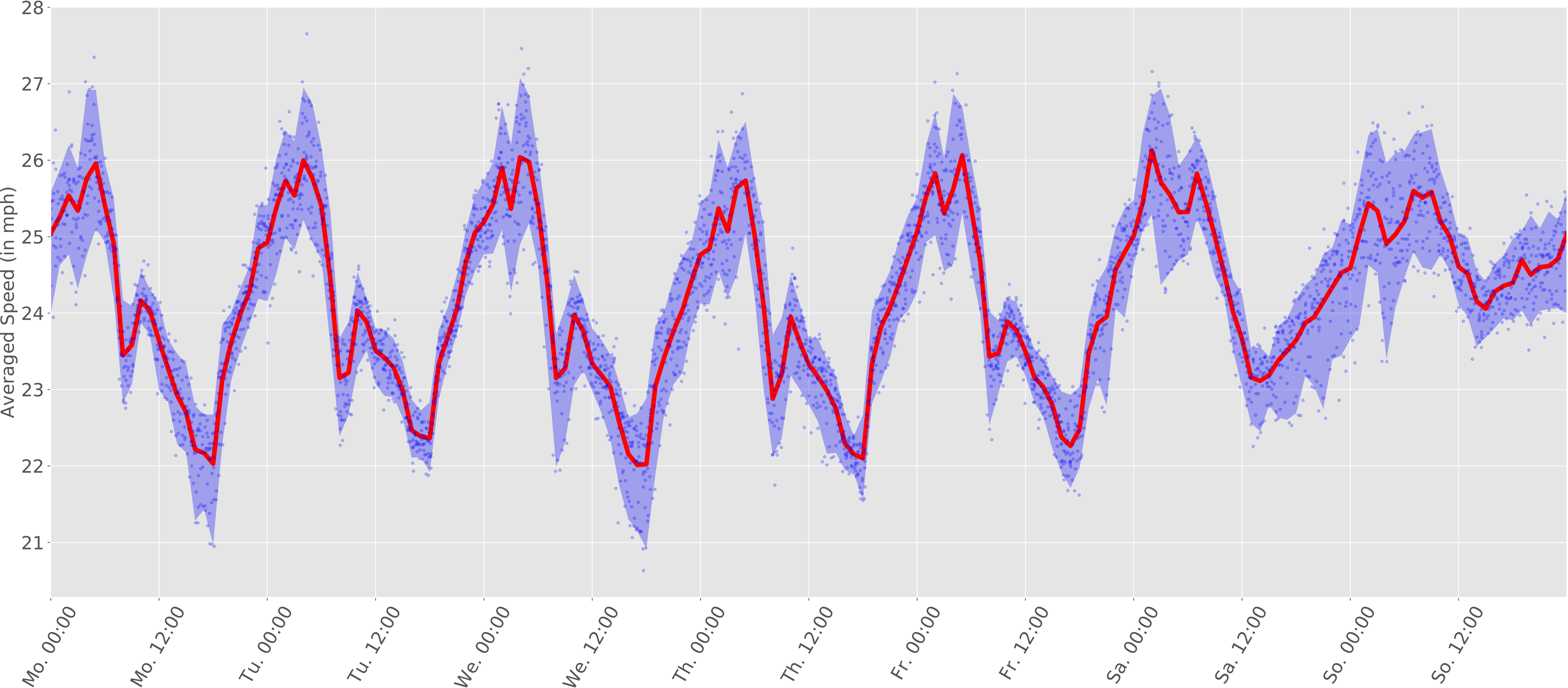}
\caption{Averaged recorded travel time data.}\label{fig:rawdata}
\end{figure}

Figure~\ref{fig:rawdata} visualizes the travel time data used in these experiments plotted against the time of one week, where each point represents the average travel time over all segments in the network in one observation. The red line shows the hourly average travel time, and the blue shaded area represents the corresponding 95\% confidence band.

As segments are purely geographical objects without structure, we needed to create a graph for our experiments. To this end, segments were split when they crossed or nearly crossed, and start- and end-points that were sufficiently close to each other were identified as the same node. The resulting graph is shown in Figure~\ref{chicago2}; note that this process slightly simplified the network, but kept its structure intact. The final graph contains 538 nodes and 1308 arcs. Using arc length and speed, we calculate their respective traversal time for each of the 4363 data points. In the following, we refer to 4363 \textit{scenarios} generated this way.

We then used this full dataset to derive the following subsets aimed at providing robust solutions in different contexts:
\begin{enumerate}
\item Find a path that is robust when driving during morning rush hours. We only use scenarios sampled on weekdays from 8am to 10am. These are 271 such scenarios (``mornings dataset'').
\item Find a path that is robust when driving during evening rush hours. We only use scenarios sampled on weekdays from 4pm to 6pm. These are 272 scenarios (``evenings dataset'').
\item Find a path that is robust when driving during a Tuesday. There are 671 scenarios sampled on Tuesdays (``Tuesdays dataset'').
\item Find a path that is robust when driving during the weekend. There are 1141 scenarios sampled on Saturdays and Sundays (``weekends dataset'').
\item Find a path that is robust when no additional information is given. We use all 4363 scenarios (``complete dataset'').
\end{enumerate}
In the following, we present results only for the mornings dataset. Results for the other datasets can be found in \ref{appendix}.

\subsection{Setup}

Each uncertainty set is equipped with a scaling parameter. For each parameter we generated 20 possible values, reflecting a reasonable range of choices for a decision maker:
\begin{itemize}
\item For $\cU^{CH}$ and $\cU^I$, $\lambda \in \{ 0.05, 0.10, \ldots, 1.00\}$.
\item For $\cU^E$, $\lambda \in \{ 0.5, 1.0, \ldots, 10.0\}$.
\item For $\cU^B$, $\Gamma \in \{ 1, 2, \ldots, 20\}$.
\item For $\cU^{PH}$, we used columns $\pmb{q}_1, \pmb{q}_3, \ldots, \pmb{q}_{39}$.
\item For $\cU^{SPH}$, we used columns $\pmb{q}_1, \pmb{q}_2, \ldots, \pmb{q}_{20}$.
\end{itemize}
Additionally, we calculate a solution to the average-case scenario $\hat{\pmb{c}}$. Note that this is a special case of all uncertainty sets presented here if the scaling parameter is sufficiently small. Each uncertainty set is generated using $75\%$ of scenarios sampled uniformly (e.g., 203 out of 271), and we evaluate solutions in-sample and out-sample separately. Furthermore, we generated 200 random $s-t$ pairs uniformly over the node set, and used each of the $6\cdot 20$ methods on the same 200 pairs. Each of our 120 methods, hence, generates $200\cdot 271 = 54,200$ objective values for the mornings set.

It is non-trivial to assess the quality of these robust solutions, see \cite{chassein2016performance}. If one just uses the average objective value, as an example, then one could as well calculate the solution optimizing the average scenario case to find the best performance with respect to this measure. To find a balanced evaluation of all methods, we used three performance criteria:
\begin{itemize}
\item the average objective value over all $s-t$ pairs and all scenarios,
\item the average of the worst-case objective value for each $s-t$ pair, and
\item the average value of the worst 5\% of objective values for each $s-t$ pair (as in the CVaR measure)
\end{itemize}
Note that many more criteria would be possible to use.

For all experiments we used a computer with a 16-core Intel Xeon E5-2670 processor, running at 2.60 GHz with 20MB cache, and Ubuntu 12.04. Processes were pinned to one core. We used Cplex v.12.6 to solve all problem formulations (note that specialized combinatorial algorithms are available for some problems).

\subsection{Results}

\begin{figure}[htbp]
\centering
\subfigure[in-sample]{\label{av-max-in}\includegraphics[width=\textwidth]{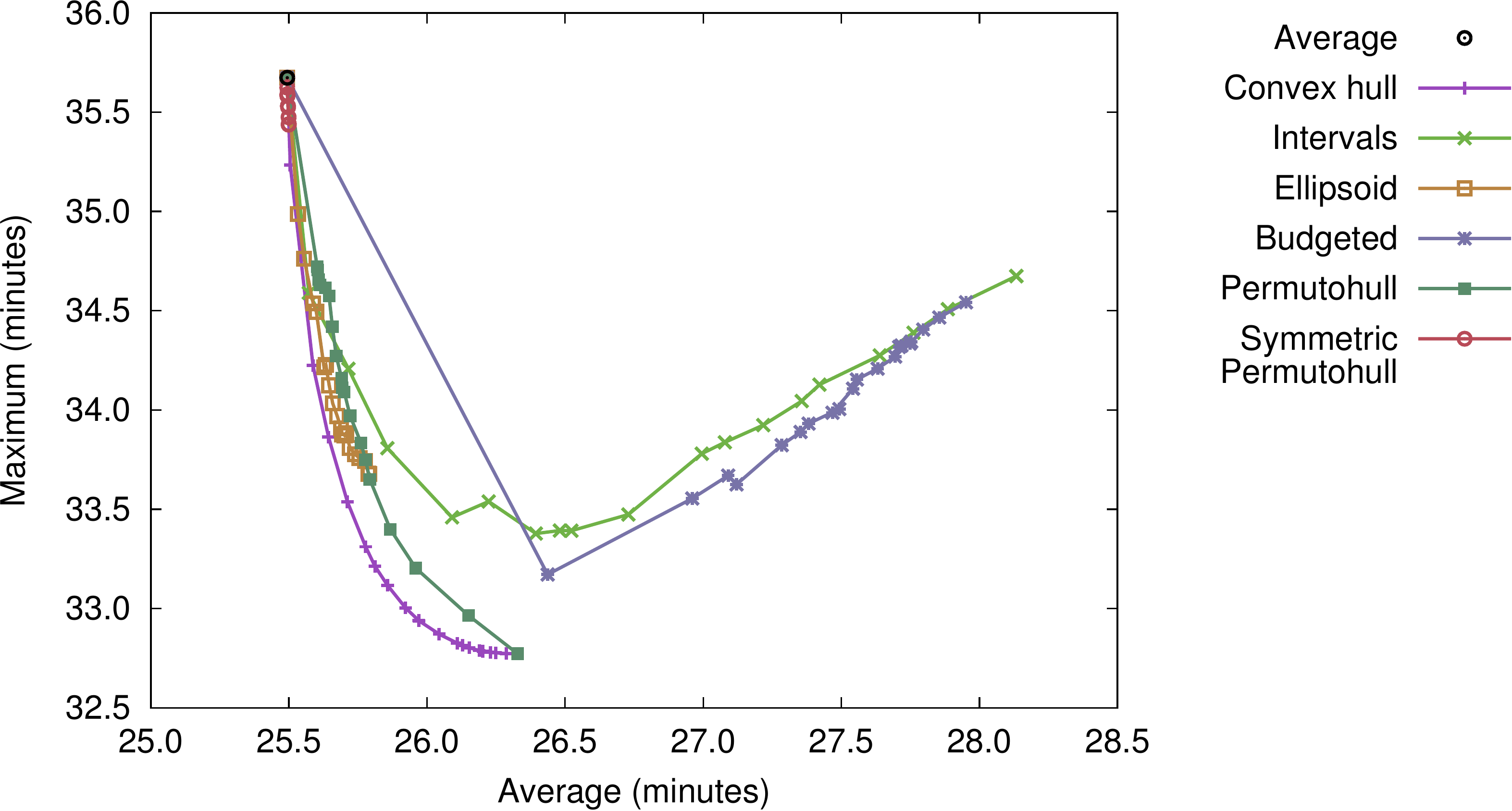}}  
    \hfill
\subfigure[out-sample]{\label{av-max-out}\includegraphics[width=\textwidth]{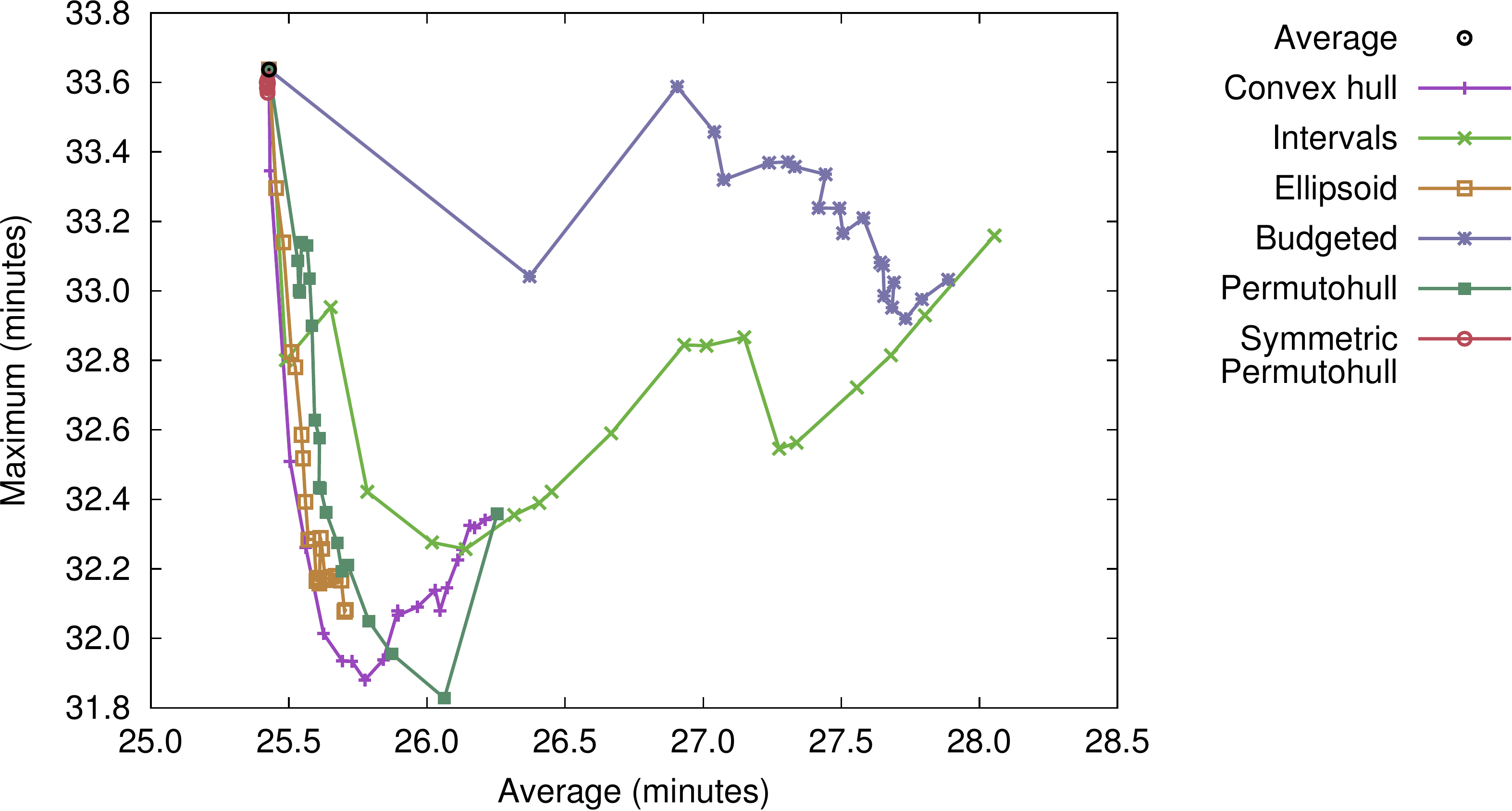}}      
\caption{Average vs worst-case performance.}\label{results-av-max}
\end{figure}

\begin{figure}[htbp]
\centering
\subfigure[in-sample]{\label{av-cvar-in}\includegraphics[width=\textwidth]{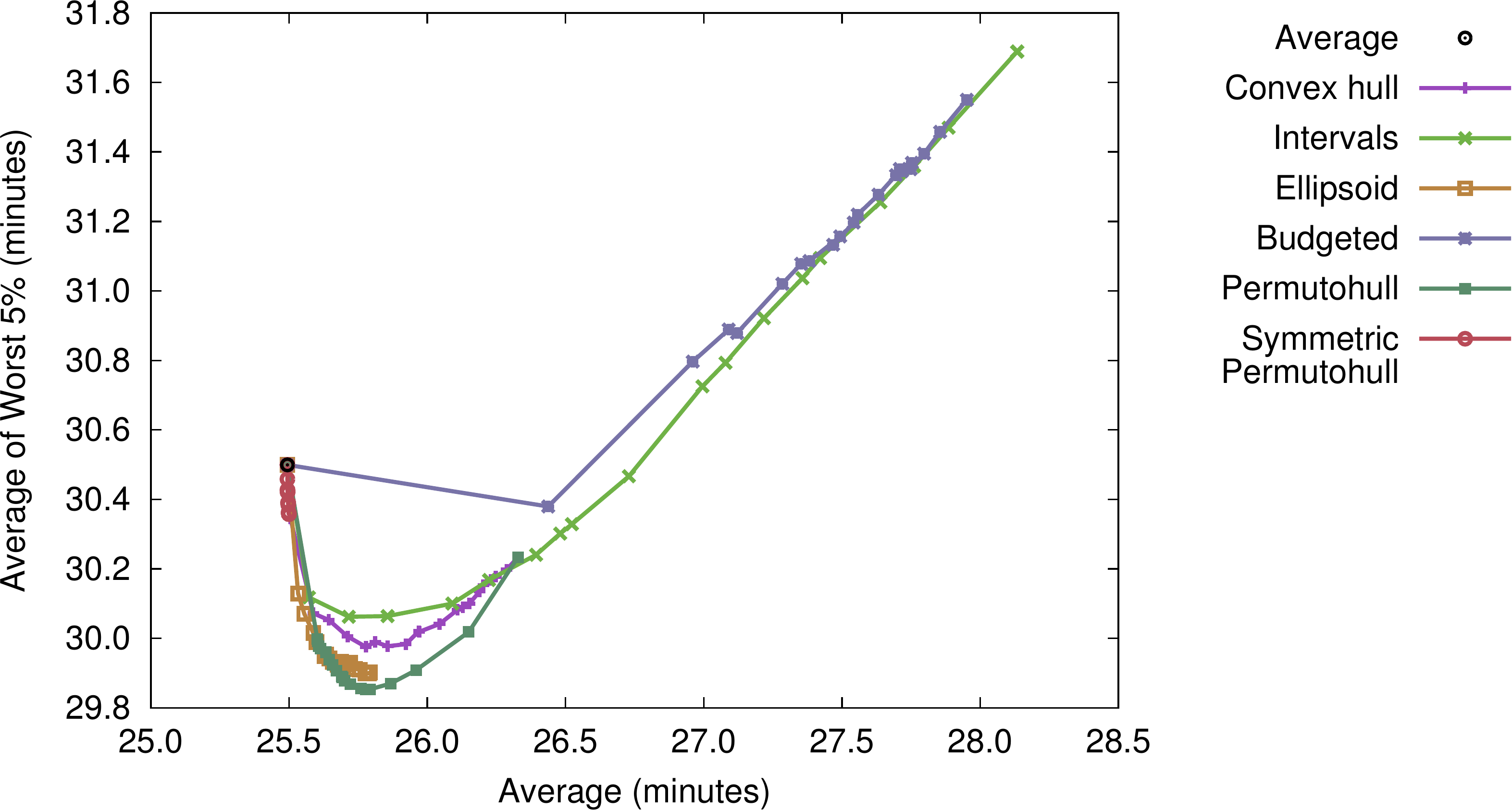}}      
    \hfill
\subfigure[out-sample]{\label{av-cvar-out}\includegraphics[width=\textwidth]{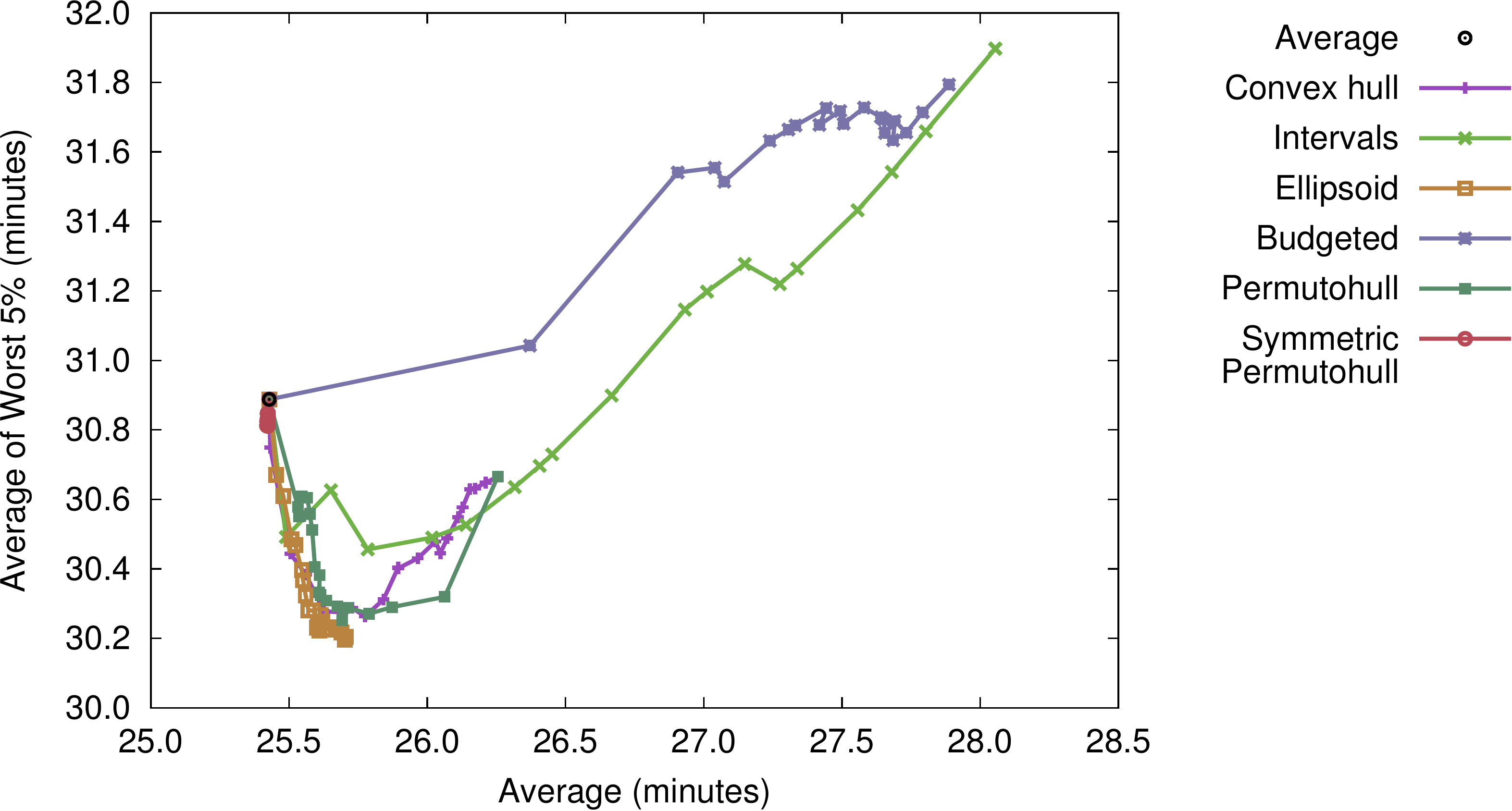}}      
\caption{Average vs CVaR performance.}\label{results-av-cvar}
\end{figure}

\begin{figure}[htbp]
\centering
\includegraphics[width=\textwidth]{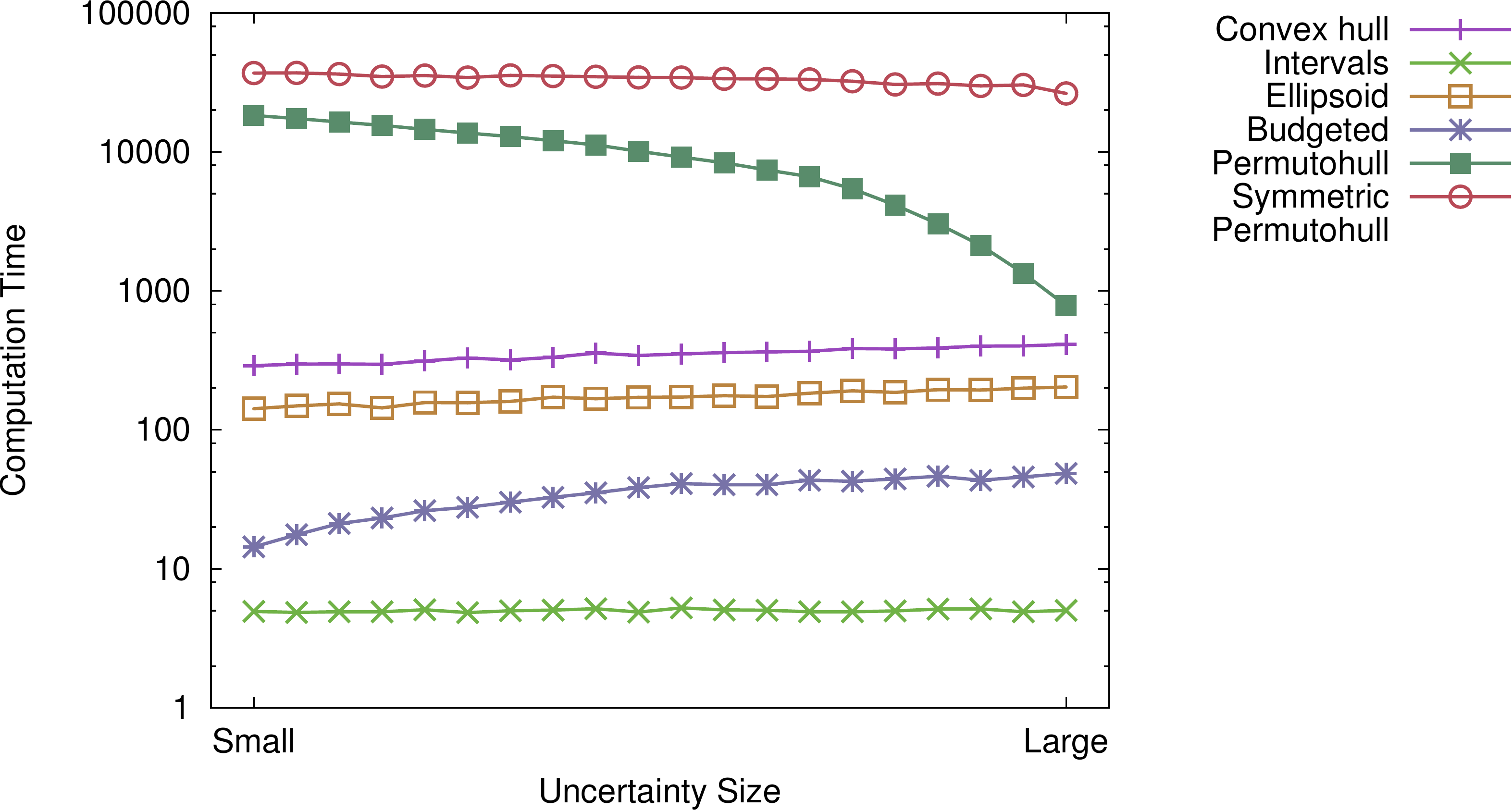}
\caption{Computation times in seconds.}\label{time}
\end{figure}

We present the performance of solutions in Figures~\ref{results-av-max} and \ref{results-av-cvar}. In each plot, the 20 parameter settings that belong to the same uncertainty set are connected by a line, including the average case as a 21st point. They are complemented with Figure~\ref{time} showing the total computation times for the methods over all 200 shortest path calculations. 

The first set of plots in Figure~\ref{results-av-max} shows the trade-off between the average and the maximum objective value; the second set of plots in Figure~\ref{results-av-cvar} shows the trade-off between the average and the average of the 5\% worst objective values. For each case, the in-sample and out-sample performance is shown. All values are in minutes of travel time. Note that for all performance measures, smaller values indicate a better performance -- hence, good trade-off solutions should move from the top left to the bottom right of the plots. In general, the points corresponding to the parameter settings that give weight to the average performance are on the left sides of the curves, while the more robust parameter settings are on the right sides, as would be expected.

We first discuss the in-sample performance in Figure~\ref{av-max-in}. In general, we find that most concepts indeed present a trade-off between average performance and robustness through their scaling parameter. Solutions calculated using the convex hull dominate the others. Symmetric permutohull solutions tend to focus on a good average performance, while ellipsoid and permutohull solutions show a broader front over the two criteria. Interval and budgeted uncertainty solutions tend to perform worse for larger scaling values, without the desired trade-off property, which confirms previous findings in \cite{chassein2016bicriteria}.

Comparing these findings with the out-sample results in Figure~\ref{av-max-out}, we see that a general ranking of concepts is kept intact. The solutions generated by the convex hull lose their trade-off property, as they are apparently over-fitted to the data used in the sample (see also the results in \ref{appendix}). We also find that the interval uncertainty outperforms budgeted uncertainty here, while they showed similar performance in-sample.

We now consider the results presented in Figure~\ref{results-av-cvar}. Here the average is plotted against the average performance of the 5\% worst performing scenarios, averaged over all $s-t$ pairs. While the convex hull solutions showed the best trade-off in Figure~\ref{av-max-in}, we find that the permutohull solutions are prominent among the non-dominated points in this case. As before, symmetric permutohull solutions tend to remain on one end of the spectrum with good average performance. Solutions based on ellipsoidal uncertainty are still among the best-performing approaches, and stable when considered out-sample (compare Figures~\ref{av-cvar-in} and \ref{av-cvar-out}). Interval uncertainty and in particular budgeted uncertainty do not perform well in comparison with the other approaches.

Regarding computation times (see Figure~\ref{time}), note that the two polynomially solvable approaches (intervals and budgeted) are also the fastest when using Cplex; these computation times can be further improved using specialized algorithms. Using ellipsoids is faster than using the convex hull, which is in turn faster than using the symmetric permutohull. For the standard permutohull, the computation times are sensitive to the uncertainty size; if the $\pmb{q}$ vector that is used in the model has only few entries, computation times are smaller. This is in line with the intuition that the problem becomes easier if fewer scenarios need to be considered.

To summarize our findings in our experiment on the robust shortest path problem with real-world data:
\begin{itemize}
\item Convex hull solutions show good in-sample performance, but are not stable when facing scenarios out of sample.

\item Interval solutions do not perform well in general, but are easy and fast to compute, which makes them a reasonable approach, in particular for smaller scalings.

\item Budgeted uncertainty does not seem an adequate choice for robust shortest path problems. Scaling interval uncertainty sets gives better results and is easier to use and to solve.

\item Ellipsoidal uncertainty solutions have good and stable overall performance and represent a large part of the non-dominated points in our results.

\item Permutohull solutions offer good trade-off solutions, whereas symmetric permutohull solutions tend to be less robust, but provide an excellent average performance. These methods also require most computational effort to find.
\end{itemize}

In the light of these findings, permutohull and ellipsoidal uncertainty tend to produce solutions with the best trade-off, while being computationally more challenging than most of the other approaches.
The algorithmic research for robust shortest path problems with such structure should therefore be studied further. Results on additional experiments leading to the same conclusions can be found in the appendix.  In the following section, we consider a variant of ellipsoidal uncertainty where correlation between arcs is ignored.

%\section{Efficient Algorithm for Axis-Parallel Ellipsoids}
\section{Ellipsoidal Uncertainty Sets}\label{sec:ellipsoids}

\subsection{From General to Axis-Parallel Ellipsoids}

Since our experiments have shown that ellipsoidal uncertainty sets are a reasonable choice, we devote this section to these sets. First, we show that changing the general ellipsoid to an axis-parallel ellipsoid by setting all non-diagonal entries of $\pmb{\Sigma}$ to $0$, has almost no effect on the found solutions. Second, we derive a specialized branch-and-bound algorithm for such axis-parallel ellipsoidal uncertainty sets which clearly outperforms the standard approach of using a generic solver. Problems with the same structure were previously considered in \cite{nikolova2009high}, where a heuristic method was proposed. 

\begin{figure}[htbp]
\centering
\includegraphics[width=\textwidth]{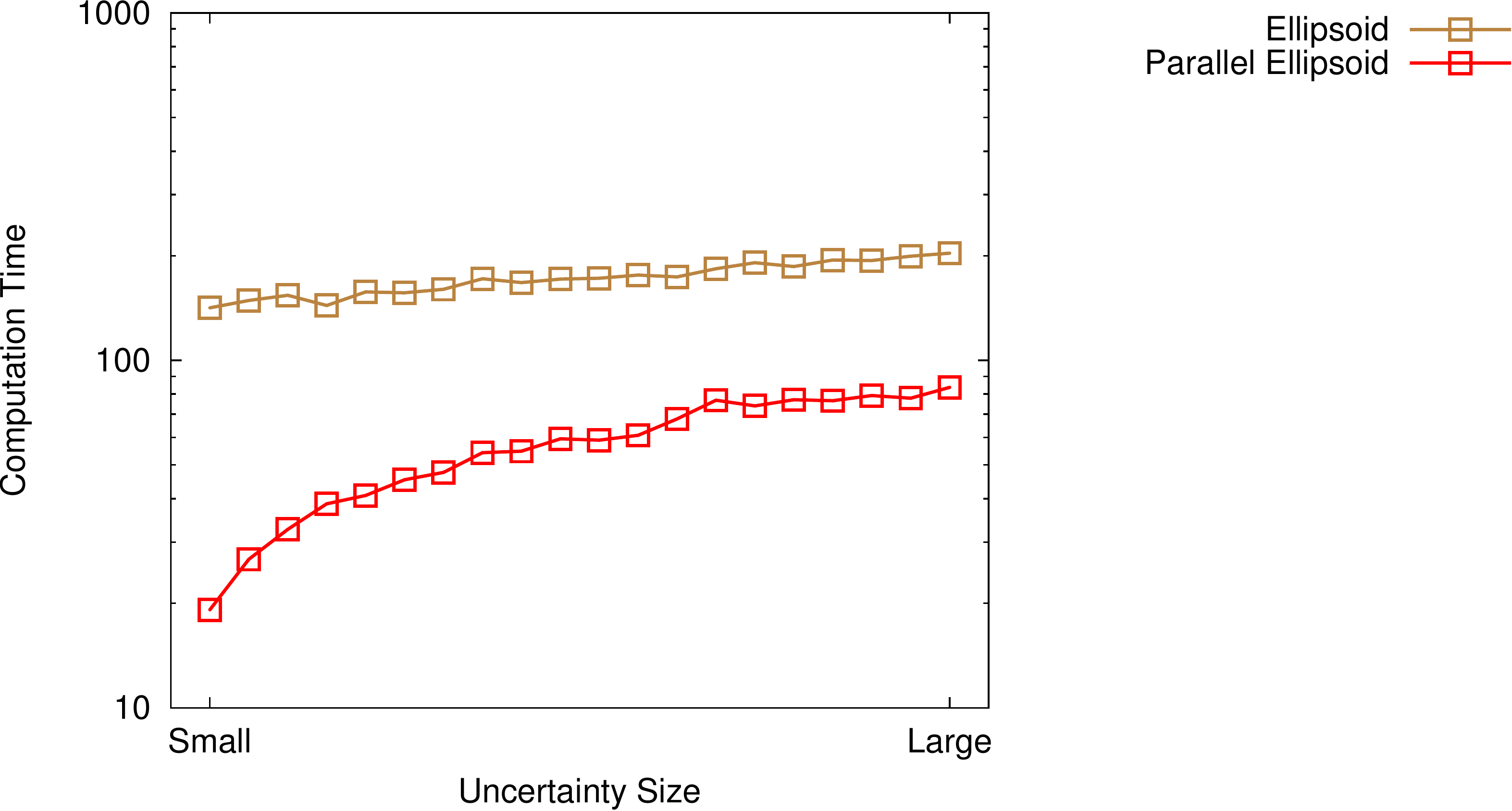}
\caption{Comparison between general and axis-parallel ellipsoids, computation times in seconds.}\label{time-2}
\end{figure}

Comparing general and axis-parallel ellipsoids for the experiments presented in Section~\ref{sec:exp}, we find that the maximum deviation over all plotted datapoints is less than $0.002\%$ for average travel times, less than $0.005\%$ for average worst-case values, and less than $0.003\%$ for average CVaR values. That is, plotted in our figures, general and axis-parallel ellipsoids would look indistinguishable. On the other hand, using axis-parallel ellipsoidal uncertainty sets in the robust model instead of general ellipsoids decreases the computation time significantly, see Figure~\ref{time-2}. The computation time can be further reduced by the use of specialized algorithms, as shown in the next section.

\subsection{Efficient Algorithm for Axis-Parallel Ellipsoids}

\subsubsection{A Bicriteria Perspective}

In this section, we describe an efficient branch-and-bound algorithm to solve the robust shortest path problem if the uncertainty set is given as an axis parallel ellipsoid. Recall that the mathematical formulation of the problem is 
\begin{align*}
\min\ & \hat{\pmb{c}}^t\pmb{x} + z \\
\text{s.t. } & z^2 \ge \left( \pmb{x}^t \pmb{\Sigma} \pmb{x}\right) \\ 
& \pmb{x}\in\X
\end{align*}
where $\pmb{\Sigma}$ is a diagonal matrix specifying the shape and the size of the ellipsoid. Since $\pmb{x}$ is a binary vector, we can simplify the quadratic expression $\pmb{x}^t \pmb{\Sigma} \pmb{x}$ to a linear expression $\pmb{d}^t\pmb{x}$, where $\pmb{d}$ is the diagonal of $\pmb{\Sigma}$. Hence the problem can be reduced to
\begin{align*}
\min\ & \hat{\pmb{c}}^t\pmb{x} + \sqrt{\pmb{d}^t\pmb{x}} \\
\text{s.t. } & \pmb{x}\in\X
\end{align*}
As pointed out in \cite{nikolova2009high}, this problem can be transformed to the following bicriteria optimization problem. 
\begin{align*}
\min\ & \begin{pmatrix}
\hat{\pmb{c}}^t\pmb{x} \\
\pmb{d}^t\pmb{x} \end{pmatrix} \\
\text{s.t. } & \pmb{x}\in\X
\end{align*}
It is shown that each optimal solution of the robust problem is an efficient extreme solution of this bicriteria optimization problem. We call a solution $\pmb{x}^*$ of this bicriteria optimization problem \emph{efficient extreme} if there exists $\alpha_0$ and $\alpha_1$ with $0\leq \alpha_0 < \alpha_1 \leq 1$ such that for all $\alpha \in [\alpha_0,\alpha_1]$ it holds that there exists no other solution $\pmb{x}'$ with $(\alpha \pmb{c} + (1-\alpha) \pmb{d})^t\pmb{x}' < (\alpha \pmb{c} + (1-\alpha) \pmb{d})^t\pmb{x}^*$. This means that we can find efficient extreme solution by solving the following weighted sum problem which corresponds to a classic shortest path problem.
\begin{align*}
\min\ & (\alpha\hat{\pmb{c}} + (1-\alpha)\pmb{d})^t\pmb{x} \\
\text{s.t. } & \pmb{x}\in\X
\end{align*}
Hence, the robust solution can be found by computing all efficient extreme solutions of the bicriteria problem. Unfortunately, there is no polynomial bound for the number of efficient extreme solutions of a bicriteria shortest path problem. In fact, it has been shown in \cite{carstensen1983complexity} that there exist instances of the bicriteria shortest path problem with a subexponential number of efficient extreme solutions. Hence, \cite{nikolova2009high} proposes a heuristic to compute only a subset of all efficient extreme solutions of the bicriteria problem. Among the found solutions, the best is chosen with respect to the robust objective function. 

In the following, we present an exact algorithm which is guaranteed to find an efficient extreme solution which is optimal for the robust problem without computing all efficient extreme solutions. Unfortunately, we cannot prove that the number of computed solutions by the exact algorithm is polynomially bounded. However, for real-world or randomly generated instances the number of computed solutions is so small that it can be assumed to be constant. We verify this claim in computational experiments. For convenience, we first present a naive algorithm to compute the complete set of all extreme efficient solutions.

\subsubsection{Naive Algorithm}
\label{sec:nai_alg}

\begin{algorithm}
\caption{Naive Algorithm to Compute all Efficient Extreme Solutions}\label{naive_alg}
\begin{algorithmic}[1]
\State Compute $\pmb{x}_l = \operatorname{arglexmin}_{\pmb{x} \in \X} (\hat{\pmb{c}}^t \pmb{x},\pmb{d}^t \pmb{x})$.
\State Compute $\pmb{x}_r = \operatorname{arglexmin}_{\pmb{x} \in \X} (\pmb{d}^t \pmb{x},\hat{\pmb{c}}^t \pmb{x})$.
\State \Return EXPLORE($\pmb{x}_l,\pmb{x}_r$)
\end{algorithmic}
\end{algorithm}

Note that the $\operatorname{lexmin}$ in Step~$1$ (or Step~$2$, respectively) of Algorithm~\ref{naive_alg} can be found by solving a problem of the form $\min_{\pmb{x} \in \X} ( (1-\epsilon) \hat{\pmb{c}} + \epsilon \pmb{d})^t \pmb{x}$ for a sufficiently small chosen $\epsilon$. The subroutine desribed as Algorithm~\ref{naive_subroutine} recursively finds all efficient extreme solutions. The recursion halts if the found solution is not efficient extreme anymore. 

\begin{algorithm}
\caption{Recursive Subroutine for the Naive Algorithm}\label{naive_subroutine}
\begin{algorithmic}[1]
\Procedure{EXPLORE}{$\pmb{x}_0,\pmb{x}_1$}
\State Set $\alpha_m$ such that $(\alpha_m \hat{\pmb{c}} + (1-\alpha_m) \pmb{d})^t \pmb{x}_0 = (\alpha_m \hat{\pmb{c}} + (1-\alpha_m) \pmb{d})^t \pmb{x}_1$.
\State Compute $\pmb{x}^* = \operatorname{argmin}_{\pmb{x} \in \X} ( \alpha_m \hat{\pmb{c}} + (1-\alpha_m) \pmb{d})^t \pmb{x}$
\If {$(\alpha_m \hat{\pmb{c}} + (1-\alpha_m) \pmb{d})^t \pmb{x}^* < (\alpha_m \hat{\pmb{c}} + (1-\alpha_m) \pmb{d})^t \pmb{x}_0$} 
\State \Return  EXPLORE($\pmb{x}_0,\pmb{x}^*$) $\cup \ \{\pmb{x}^*\} \ \cup$ EXPLORE($\pmb{x}^*,\pmb{x}_1$).
\Else
\State \Return $\emptyset$
\EndIf
\EndProcedure
\end{algorithmic}
\end{algorithm}

We remark that for each efficient extreme solution $\pmb{x}^*$, it is guaranteed that Algorithm~\ref{naive_alg} either finds $\pmb{x}^*$ or an alternative solution $\pmb{x}_{\text{alt}}$ with $\hat{\pmb{c}}^t\pmb{x}_{\text{alt}} = \hat{\pmb{c}}^t\pmb{x}^*$ and $\pmb{d}^t\pmb{x}_{\text{alt}} = \pmb{d}^t\pmb{x}^*$. Further, it is not guaranteed that all solutions returned by Algorithm~\ref{naive_alg} are efficient extreme. However, all non efficient extreme solutions could be easily removed. 

\subsubsection{Improved Algorithm}
\label{sec:imp_alg}

For the improved algorithm, we first find the two lexicographic minimal solutions $\pmb{x}_l$ and $\pmb{x}_r$ as in the naive method (see Figure~\ref{alg_fig_1}). We denote by $p: \X \rightarrow \R^2, p(\pmb{x}) = (p_1(\pmb{x}),p_2(\pmb{x})) =  \left(\hat{\pmb{c}}^t\pmb{x},\pmb{d}^t\pmb{x}\right)$ the map from the solution space to the two-dimensional objective space of the bicriteria optimization problem. 

%The idea of the improved algorithm is to maintain a list of intervals $\mathcal{L}$ such that all efficient extreme solutions $\pmb{x}^*$ which are not found yet and which may have a better objective value (with respect to the robust objective function) than the currently best solution fulfill that $p(\pmb{x}^*)_1 \in I$ where $I \in \mathcal{L}$.  To initialize $\mathcal{L}$ we make the following considerations.

From the definition of efficient extreme solutions, it follows that 
$p(\pmb{x}^*)$ is contained in the triangle $\Delta_{\text{unexp}}$ with the vertices $p(\pmb{x}_l), (p_1(\pmb{x}_l),p_2(\pmb{x}_r)), $ and  $p(\pmb{x}_r)$
for all efficient extreme solutions $\pmb{x}^*$ (see Figure~\ref{alg_fig_1}). Denote by $\pmb{x}_{\text{rob}}^* = \operatorname{argmin} (\hat{\pmb{c}}^t\pmb{x}_l + \sqrt{\pmb{d}^t\pmb{x}_l},\hat{\pmb{c}}^t\pmb{x}_r + \sqrt{\pmb{d}^t\pmb{x}_r})$ the current best solution for the robust problem and by $OBJ$ the corresponding objective value. Note that for all solutions $\pmb{x}$ which could improve the actual best solution, it must hold that $p_2(\pmb{x}) < (OBJ - p_1(\pmb{x}))^2$, i.e., $p(\pmb{x})$ must be contained the in the region $\mathcal{R}_{\text{imp}} = \{z \in \R^2 \mid z_2 < (OBJ-z_1)^2 \}$ (see Figure~\ref{alg_fig_2}).

\begin{figure}[htbp]
\centering
\subfigure[To initialize the algorithm we compute the two solutions which minimize the first and second objective function. The unexplored region which might contain efficient extreme solutions is marked with red diagonal lines.]{\hspace{0.04\textwidth}\label{alg_fig_1} \begin{tikzpicture}[scale=0.7,domain=0:4]
\draw[-,color=white] (0.5,-0.2) -- (0.5,0.2);
\node at (0.5,4) [circle,fill=black,inner sep = 0pt, minimum size=4pt] {};
\node at (3,0.5) [circle,fill=black,inner sep = 0pt, minimum size=4pt] {};
\fill [pattern=north east lines, pattern color = red] (0.5,0.5)--(3,0.5)--(0.5,4)--cycle;
\draw[-,color=black] (0.5,0.5) -- (3,0.5);
\draw[-,color=black] (3,0.5) -- (0.5,4);
\draw[-,color=black] (0.5,4) -- (0.5,0.5);
\draw[-,color=black] (0,0.5) -- (4,0.5);
\draw[-,color=black] (0.5,0) -- (0.5,6.5);
\draw[->,color=black] (0.5,4) -- (0.1,4);
\draw[->,color=black] (3,0.5) -- (3,0.1);
\draw[->,thin,color=black] (0,0) -- (4,0);
\draw[->,thin,color=black] (0,0) -- (0,6.5);
\end{tikzpicture}\hspace{0.04\textwidth}} \  
\subfigure[The region of improvement which might contain solutions that improve the actual best solution is marked with blue diagonal lines. The intersection of both regions, shown in gray, defines the area in which we try to find solutions.]{\hspace{0.04\textwidth}\label{alg_fig_2}\begin{tikzpicture}[scale=0.7,domain=0:4]
\draw[-,color=white] (0.5,-0.2) -- (0.5,0.2);
\begin{scope}
\path [clip] (0.5,0.5)--(3,0.5)--(0.5,4)--cycle;
\fill [gray, domain=0:2.5, variable=\x] (0, 0)-- plot ({\x}, {(2.5-\x)*(2.5-\x)})-- (3, 0)-- cycle;
\end{scope}
\fill [pattern=north west lines, pattern color = blue , domain=0:2.5, variable=\x] (0, 0)-- plot ({\x}, {(2.5-\x)*(2.5-\x)})-- (3, 0)-- cycle;
\fill [pattern=north east lines, pattern color = red] (0.5,0.5)--(3,0.5)--(0.5,4)--cycle;
\draw[-,color=black] (0.5,0.5) -- (3,0.5);
\draw[-,color=black] (3,0.5) -- (0.5,4);
\draw[-,color=black] (0.5,4) -- (0.5,0.5);
\draw [domain=0:2.5, variable=\x] plot ({\x}, {(2.5-\x)*(2.5-\x)});
\node at (0.5,4) [circle,fill=black,inner sep = 0pt, minimum size=4pt] {};
\node at (3,0.5) [circle,fill=black,inner sep = 0pt, minimum size=4pt] {};
\draw[->,thin,color=black] (0,0) -- (4,0);
\draw[->,thin,color=black] (0,0) -- (0,6.5);
\end{tikzpicture}\hspace{0.04\textwidth}} \
\subfigure[Projecting the gray area to the first axis defines the first interval the algorithm tries to shrink. The algorithm computes the midpoint of the interval and projects it to the parabola to obtain the next direction of optimization]{\hspace{0.04\textwidth}\label{alg_fig_3}\begin{tikzpicture}[scale=0.7,domain=0:4]
\begin{scope}
\path [clip] (0.5,0.5)--(3,0.5)--(0.5,4)--cycle;
\fill [gray, domain=0:2.5, variable=\x] (0, 0)-- plot ({\x}, {(2.5-\x)*(2.5-\x)})-- (3, 0)-- cycle;
\end{scope}
\draw[->,thin,color=black] (0,0) -- (4,0);
\draw[->,thin,color=black] (0,0) -- (0,6.5);
\draw[-,color=black] (0.5,-0.2) -- (0.5,0.2);
\draw[-,color=black] (1.793,-0.2) -- (1.793,0.2);
\draw[-,dotted,color=black] (1.1464,0) -- (1.1464,1.832);
\draw[-,color=black] (0,4.9354) -- (1.8231,0);
\draw[->,color=black] (1.1464,1.832) -- (0.33427,1.532);
\node at (0.5,4) [circle,fill=black,inner sep = 0pt, minimum size=4pt] {};
\node at (3,0.5) [circle,fill=black,inner sep = 0pt, minimum size=4pt] {};
\end{tikzpicture}\hspace{0.04\textwidth}}    
  
\subfigure[A new solution is found. This solution intersects the previous unexplored region with one half-space, which results in two smaller triangles.]{\hspace{0.04\textwidth}\label{alg_fig_4}\begin{tikzpicture}[scale=0.7,domain=0:4]
\draw[-,color=white] (0.5,-0.2) -- (0.5,0.2);
\fill [pattern=north east lines, pattern color = red] (1.1464,2)--(1.7,0.5)--(3,0.5)--cycle;
\fill [pattern=north east lines, pattern color = red] (1.1464,2)--(0.5,3.74985)--(0.5,4)--cycle;
\draw[-,color=black] (1.1464,2) -- (3,0.5);
\draw[-,color=black] (1.7,0.5) -- (3,0.5);
\draw[-,color=black] (1.1464,2) -- (0.5,4);
\draw[-,color=black] (0.5,3.74985) -- (0.5,4);
\draw[->,thin,color=black] (0,0) -- (4,0);
\draw[->,thin,color=black] (0,0) -- (0,6.5);
\node at (0.5,4) [circle,fill=black,inner sep = 0pt, minimum size=4pt] {};
\node at (3,0.5) [circle,fill=black,inner sep = 0pt, minimum size=4pt] {};
\node at (1.1464,2) [circle,fill=black,inner sep = 0pt, minimum size=4pt] {};
\draw[-,color=black] (0,5.10341944) -- (1.8852,0);
\draw[->,color=black] (1.1464,2) -- (0.33427,1.7);
\end{tikzpicture}\hspace{0.04\textwidth}} \
\subfigure[Intersecting the region of improvement with the new unexplored region leads to two small areas.]{\hspace{0.04\textwidth}\label{alg_fig_5}\begin{tikzpicture}[scale=0.7,domain=0:4]
\draw[-,color=white] (0.5,-0.2) -- (0.5,0.2);
\begin{scope}
\path [clip] (1.1464,2)--(1.7,0.5)--(3,0.5)--cycle;
\fill [gray, domain=0:2.5, variable=\x] (0, 0)-- plot ({\x}, {(2.5-\x)*(2.5-\x)})-- (3, 0)-- cycle;
\end{scope}
\begin{scope}
\path [clip] (1.1464,2)--(0.5,3.74985)--(0.5,4)--cycle;
\fill [gray, domain=0:2.5, variable=\x] (0, 0)-- plot ({\x}, {(2.5-\x)*(2.5-\x)})-- (3, 0)-- cycle;
\end{scope}
\fill [pattern=north west lines, pattern color = blue , domain=0:2.5, variable=\x] (0, 0)-- plot ({\x}, {(2.5-\x)*(2.5-\x)})-- (3, 0)-- cycle;
\fill [pattern=north east lines, pattern color = red] (1.1464,2)--(1.7,0.5)--(3,0.5)--cycle;
\fill [pattern=north east lines, pattern color = red] (1.1464,2)--(0.5,3.74985)--(0.5,4)--cycle;
\draw[-,color=black] (1.1464,2) -- (3,0.5);
\draw[-,color=black] (1.1464,2) -- (1.7,0.5);
\draw[-,color=black] (1.7,0.5) -- (3,0.5);
\draw[-,color=black] (1.1464,2) -- (0.5,4);
\draw[-,color=black] (1.1464,2) -- (0.5,3.74985);
\draw[-,color=black] (0.5,3.74985) -- (0.5,4);
\draw [domain=0:2.5, variable=\x] plot ({\x}, {(2.5-\x)*(2.5-\x)});
\draw[->,thin,color=black] (0,0) -- (4,0);
\draw[->,thin,color=black] (0,0) -- (0,6.5);
\node at (0.5,4) [circle,fill=black,inner sep = 0pt, minimum size=4pt] {};
\node at (3,0.5) [circle,fill=black,inner sep = 0pt, minimum size=4pt] {};
\node at (1.1464,2) [circle,fill=black,inner sep = 0pt, minimum size=4pt] {};
\end{tikzpicture}\hspace{0.04\textwidth}} \
\subfigure[Projecting the two gray areas to the first axis defines two smaller intervals. The algorithm proceeds by shrinking, splitting or removing these intervals until the gray area corresponds to the empty set]{\hspace{0.04\textwidth}\label{alg_fig_6}\begin{tikzpicture}[scale=0.7,domain=0:4]
\begin{scope}
\path [clip] (1.1464,2)--(1.7,0.5)--(3,0.5)--cycle;
\fill [gray, domain=0:2.5, variable=\x] (0, 0)-- plot ({\x}, {(2.5-\x)*(2.5-\x)})-- (3, 0)-- cycle;
\end{scope}
\begin{scope}
\path [clip] (1.1464,2)--(0.5,3.74985)--(0.5,4)--cycle;
\fill [gray, domain=0:2.5, variable=\x] (0, 0)-- plot ({\x}, {(2.5-\x)*(2.5-\x)})-- (3, 0)-- cycle;
\end{scope}
\draw[->,thin,color=black] (0,0) -- (4,0);
\draw[->,thin,color=black] (0,0) -- (0,6.5);
\node at (0.5,4) [circle,fill=black,inner sep = 0pt, minimum size=4pt] {};
\node at (3,0.5) [circle,fill=black,inner sep = 0pt, minimum size=4pt] {};
\node at (1.1464,2) [circle,fill=black,inner sep = 0pt, minimum size=4pt] {};

\draw[-,color=black] (0.5,-0.2) -- (0.5,0.2);
\draw[-,color=black] (0.73688,-0.2) -- (0.73688,0.2);
\draw[-,color=black] (1.55602,-0.2) -- (1.55602,0.2);
\draw[-,color=black] (1.793,-0.2) -- (1.793,0.2);
\end{tikzpicture}\hspace{0.04\textwidth}}
\caption{Visualization of the improved algorithm.}\label{alg_fig}
\end{figure}

Intuitively, we always have two regions in which we are interested during the algorithm. First, the unexplored region, which may contain efficient extreme solutions which we have not found yet. At the beginning this region corresponds to $\Delta_{\text{unexp}}$. Second, the improving region, corresponding to $\mathcal{R}_{\text{imp}}$,  which could contain solutions that improve our current best solution. We intersect these two regions and project the so obtained set to the first axis. This gives a list of intervals $\mathcal{L}$. The idea of the improved algorithm is then to shrink or to remove intervals from $\mathcal{L}$ until $\mathcal{L}$ is empty.

At the beginning, we intersect the triangle $\Delta_{\text{unexp}}$ and $\mathcal{R}_{\text{imp}}$ and project the area obtained this way to the first axis. This results in the interval which we use to initialize $\mathcal{L}$ (see Figure~\ref{alg_fig_3}). The idea of the improved algorithm is then to shrink, split or remove intervals from $\mathcal{L}$ until $\mathcal{L}$ is empty. 

To do so, the improved algorithm picks an interval $I=[a,b]$ from $\mathcal{L}$, computes its midpoint $m=0.5(a+b)$ and defines a local approximation of the boundary of $\mathcal{R}_{\text{imp}}$ at $(m,(OBJ-m)^2)$ (see Figure~\ref{alg_fig_3}). The slope of the obtained line $l$ is $2(m-OBJ)$. We set $\alpha_m = 1/(1+2(OBJ-m))$ to optimize in the direction which is perpendicular to $l$ (see Figure~\ref{alg_fig_3}). 

We then compute $\pmb{x}_{\text{new}} = \operatorname{argmin}_{\pmb{x} \in \X} (\alpha_m \hat{\pmb{c}} + (1-\alpha_m) \pmb{d})^t \pmb{x}$. After we have found $\pmb{x}_{\text{new}}$ we know that for all other efficient extreme solutions $\pmb{x}^*$ it must hold that $p(\pmb{x}^*)$ must lie above the line trough $p(\pmb{x}_{\text{new}})$ with slope $2(m-OBJ)$. Hence, we can exclude a half space from the unexplored region (see Figure~\ref{alg_fig_4}) and shrink, split or remove intervals contained in $\mathcal{L}$ (see Figure~\ref{alg_fig_5}). Further, it might happen that $\pmb{x}_{\text{new}}$ improves the actual best solution in this case we update $OBJ$ and $\mathcal{R}_{\text{imp}}$ and consequently all intervals in $\mathcal{L}$ .

If the algorithm has reduced $\mathcal{L}$ to the empty set the current best solution is indeed the optimal solution of the robust optimization problem.

\subsubsection{Computational Experiments for the Improved Algorithm} \label{sec:exp_imp}

\begin{figure}[htbp]
\centering
\includegraphics[width=\textwidth]{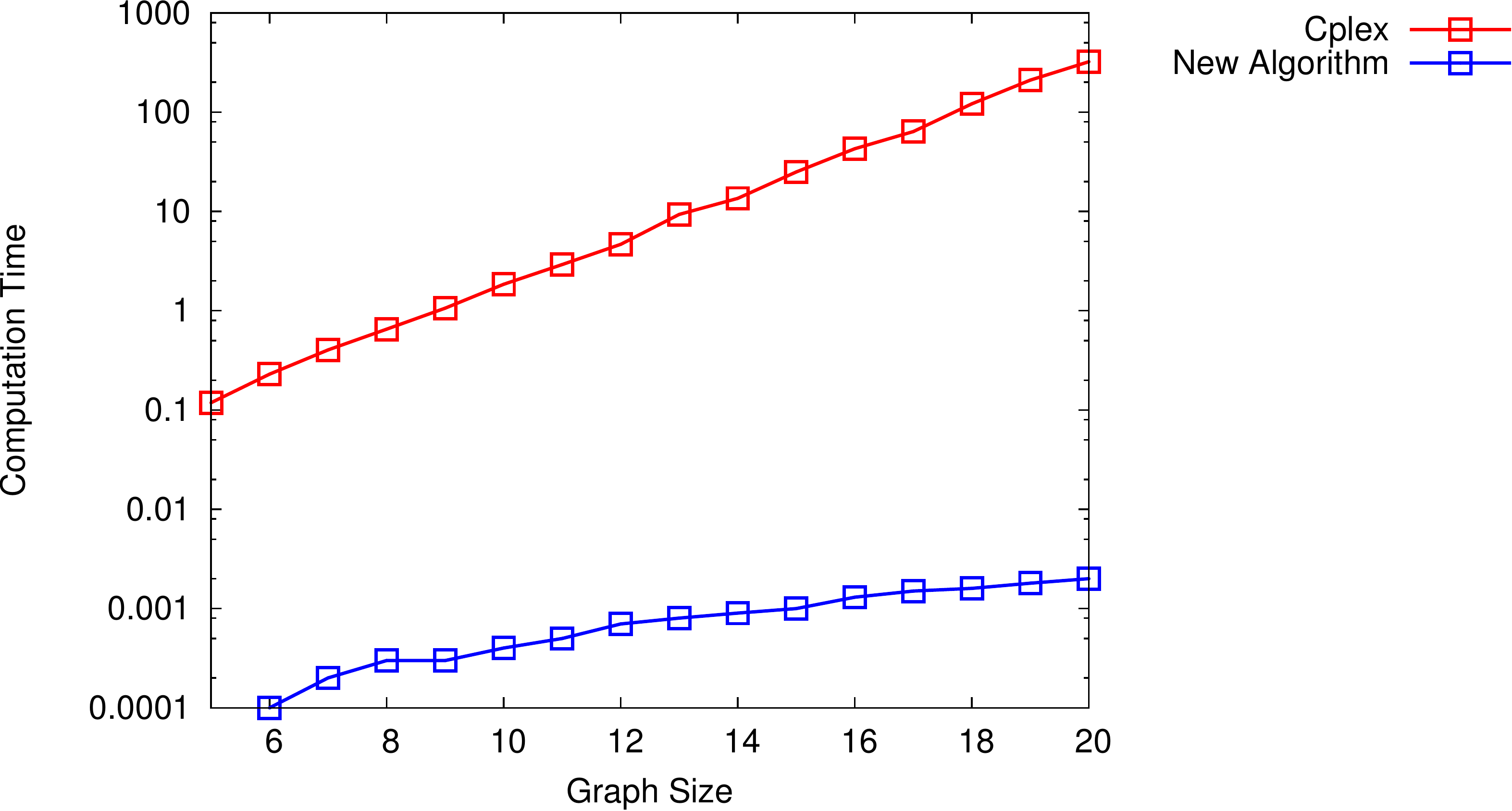}
\caption{Comparison of average computation times for grid graphs.}\label{time-3}
\end{figure}

We test the performance of the branch-and-bound algorithm on grid graphs, using the same computational environment as in Section~\ref{sec:exp}. We used the LEMON graph library (v.1.3.1) to solve the classic shortest path problems that needs to be solved during the branch-and-bound algorithm. The goal is to find a path from the upper left corner to the lower right corner of the grid. For each arc we chose $50$ values uniform at random from $[100]$. Further, we fit an axis-parallel ellipsoidal uncertainty set to these points as described previously. We vary the grid size from a $5\times 5$ grid to a $20 \times 20$ grid. For each grid size we create $100$ instances and solve them in two ways: By using Cplex to solve the resulting MISOCP and by the proposed branch-and-bound algorithm.  The averaged computation times are shown in Figure~\ref{time-3}.

It can be seen that our approach outperforms Cplex by several orders of magnitude (note the logarithmic vertical scale). While the computation times for Cplex scale exponentially with the graph size, this is not observed for our method. We show the average number of shortest path calculations required by our method in Table~\ref{tab2}, where we increase the grid size to $100\times 100$. It can be seen that on average only very few calls are required, and the increase is slow. By using further improved algorithms for shortest path calculation in road networks (see \cite{bast2016route}), the application of our method in real-time route planning is within reach.

\begin{table}[htb]
\centering\begin{tabular}{c|r}
Instance size & SP comp.\\
\hline
$5\times 5$ & 3.625 \\
$10\times 10$ & 3.929 \\
$20\times 20$ & 4.251 \\
$50\times 50$ & 4.778 \\
$100\times 100$ & 5.127
\end{tabular}
\caption{Average number of shortest path computations for different grid sizes.}\label{tab2}
\end{table}

\begin{figure}[htbp]
\centering
\includegraphics[width=\textwidth]{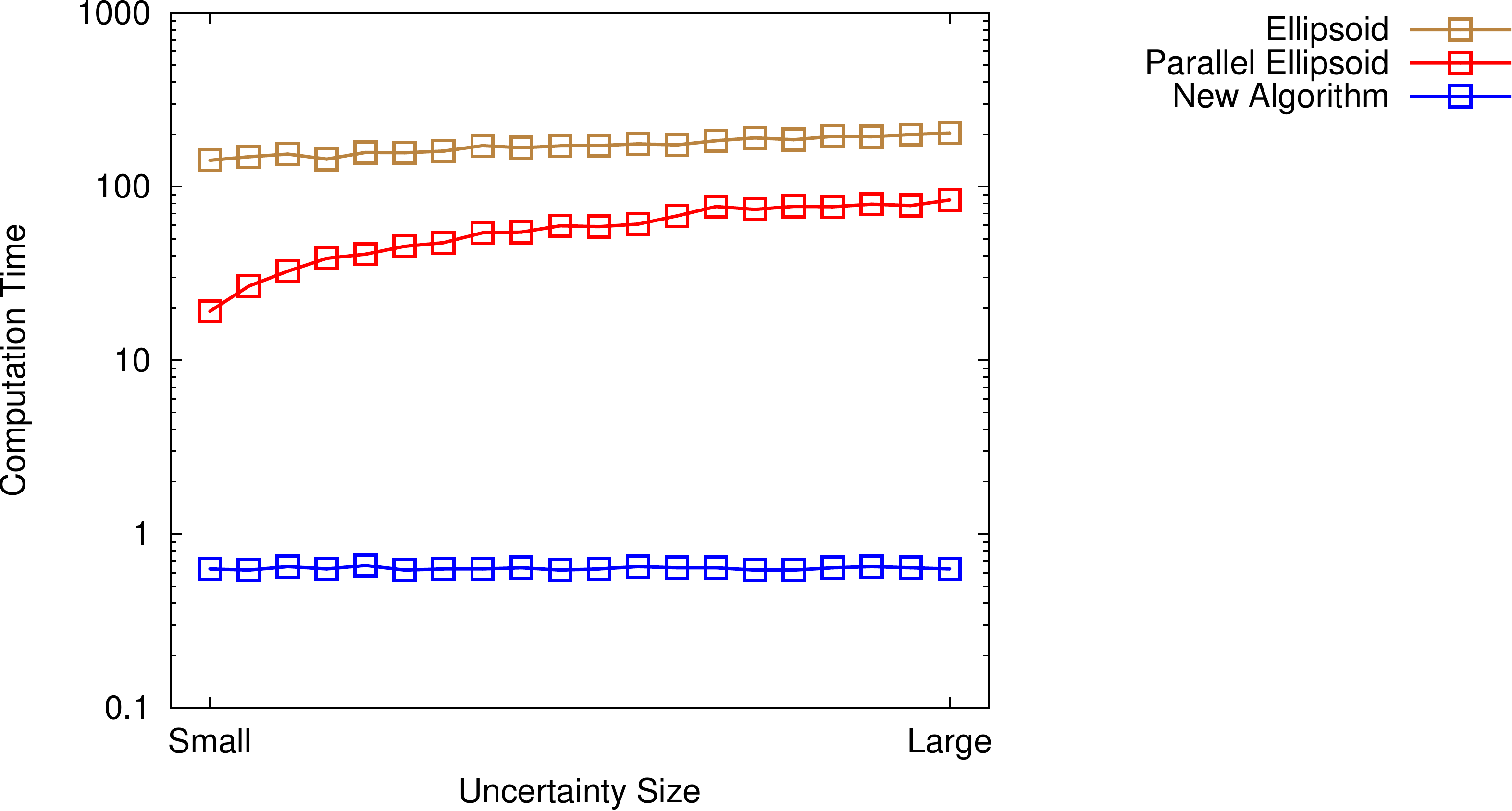}
\caption{Comparison between general and axis-parallel ellipsoids, computation times in seconds.}\label{time-4}
\end{figure}

Finally, we revisit the computation times for the real-world instance from our previous experiment (Figure~\ref{time-2}). Figure~\ref{time-4} shows the performance of our method for comparison.

\section{Conclusions}
\label{sec:conclusion} 
 
In this paper, we constructed uncertainty sets for the robust shortest path problem using real-world traffic observations for the City of Chicago. We evaluated the model suitability of these sets by finding the resulting robust paths, and comparing their in-sample and out-sample performance using different performance indicators.
Naturally, conclusions can only be drawn within the reach of the available data. It remains to be seen how the considered uncertainty sets perform on other datasets for robust shortest paths.
 
We have observed that using ellipsoidal uncertainty sets provides high-quality solutions with less computational effort than for the permutohull. If one uses only the diagonal entries of the matrix $\pmb{\Sigma}$, then one ignores the data correlation in the network, but the solution quality remains roughly the same. For the resulting problem, a specialized branch-and-bound algorithm was developed that is able to reduce computation times considerably compared to Cplex. In fact, the computational effort to solve this problem is comparable to the complexity of solving a few classic shortest path problems, which even makes the application on real-time navigation devices a possibility.

% For the resulting problem specialized algorithms exist, see, e.g. \cite{nikolova2009high}. In additional experiments we found that even by using Cplex, computation times were considerably reduced when only using the diagonal entries of $\pmb{\Sigma}$, but the solution quality remained roughly the same.

%%
%% Bibliography
%%

%% Either use bibtex (recommended), 

% \section*{References}

%% .. or use the thebibliography environment explicitely

\appendix

\section{Additional experimental results}
\label{appendix}

\begin{figure}[htbp]
\centering
\subfigure[in-sample]{\label{av-max-in-evenings}\includegraphics[width=\textwidth]{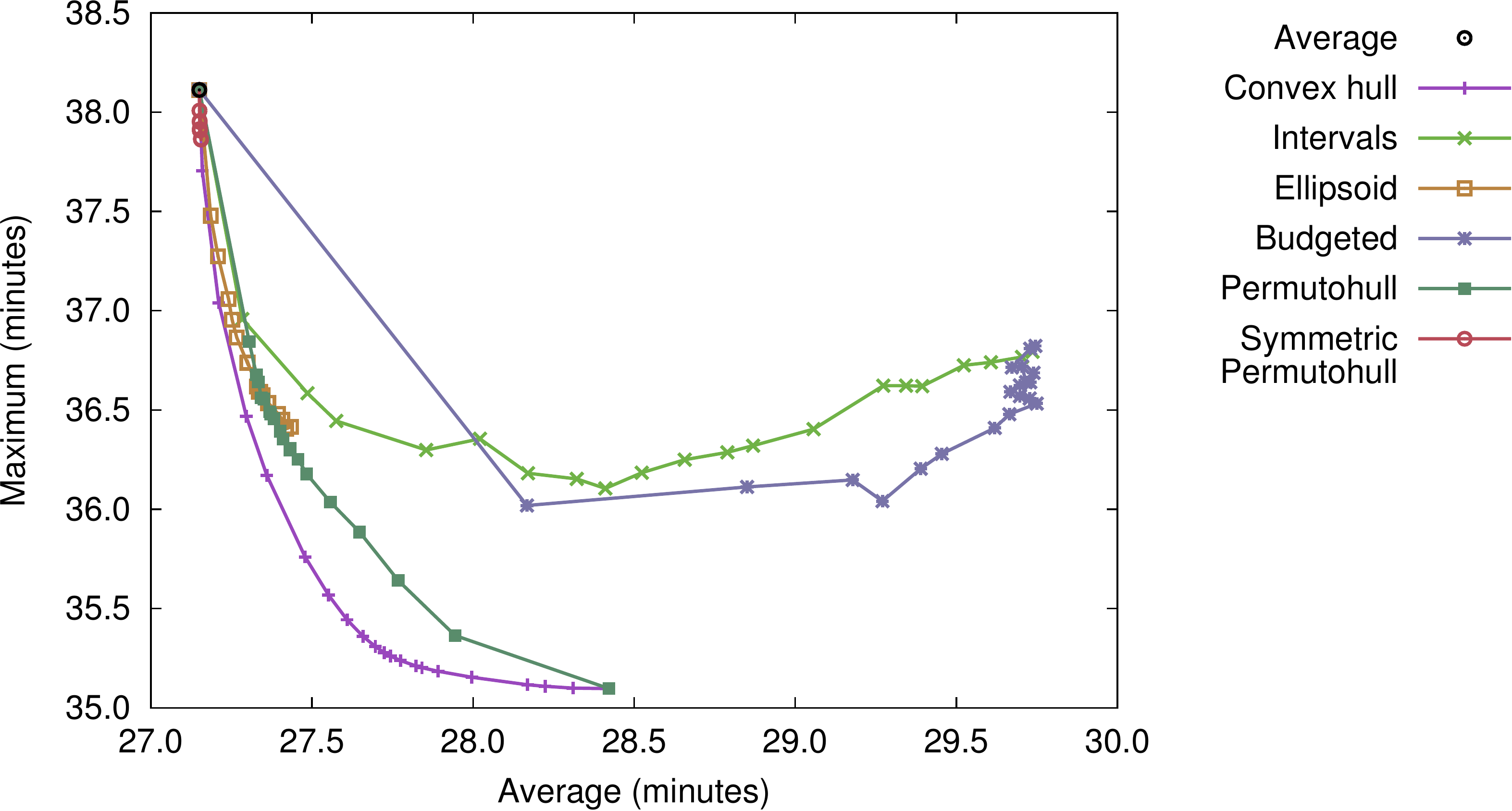}}  
    \hfill
\subfigure[out-sample]{\label{av-max-out-evenings}\includegraphics[width=\textwidth]{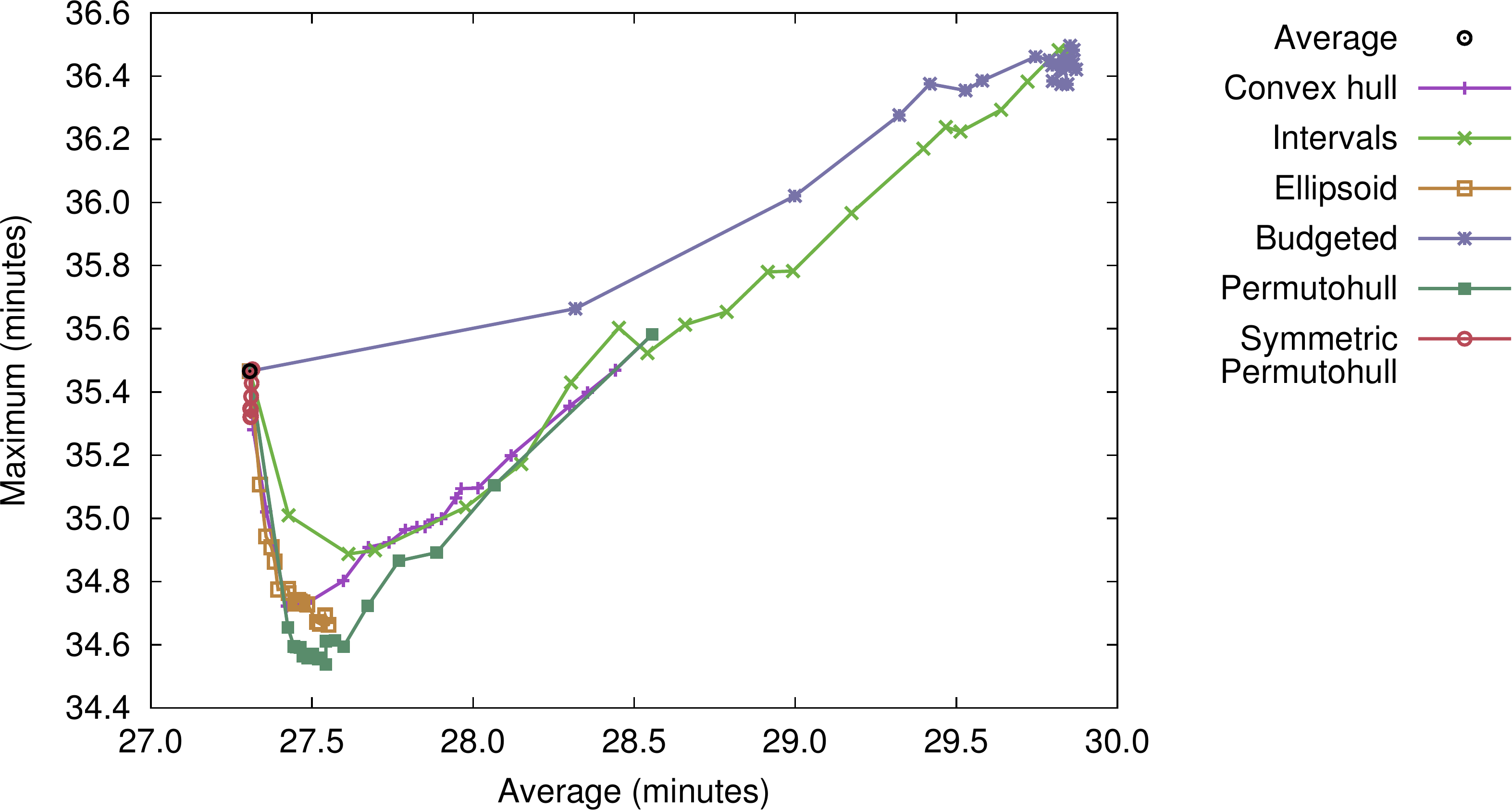}}      
\caption{Average vs worst-case performance, evenings dataset.}\label{av-max-evenings}
\end{figure}

\begin{figure}[htbp]
\centering
\subfigure[in-sample]{\label{av-cav-in-evenings}\includegraphics[width=\textwidth]{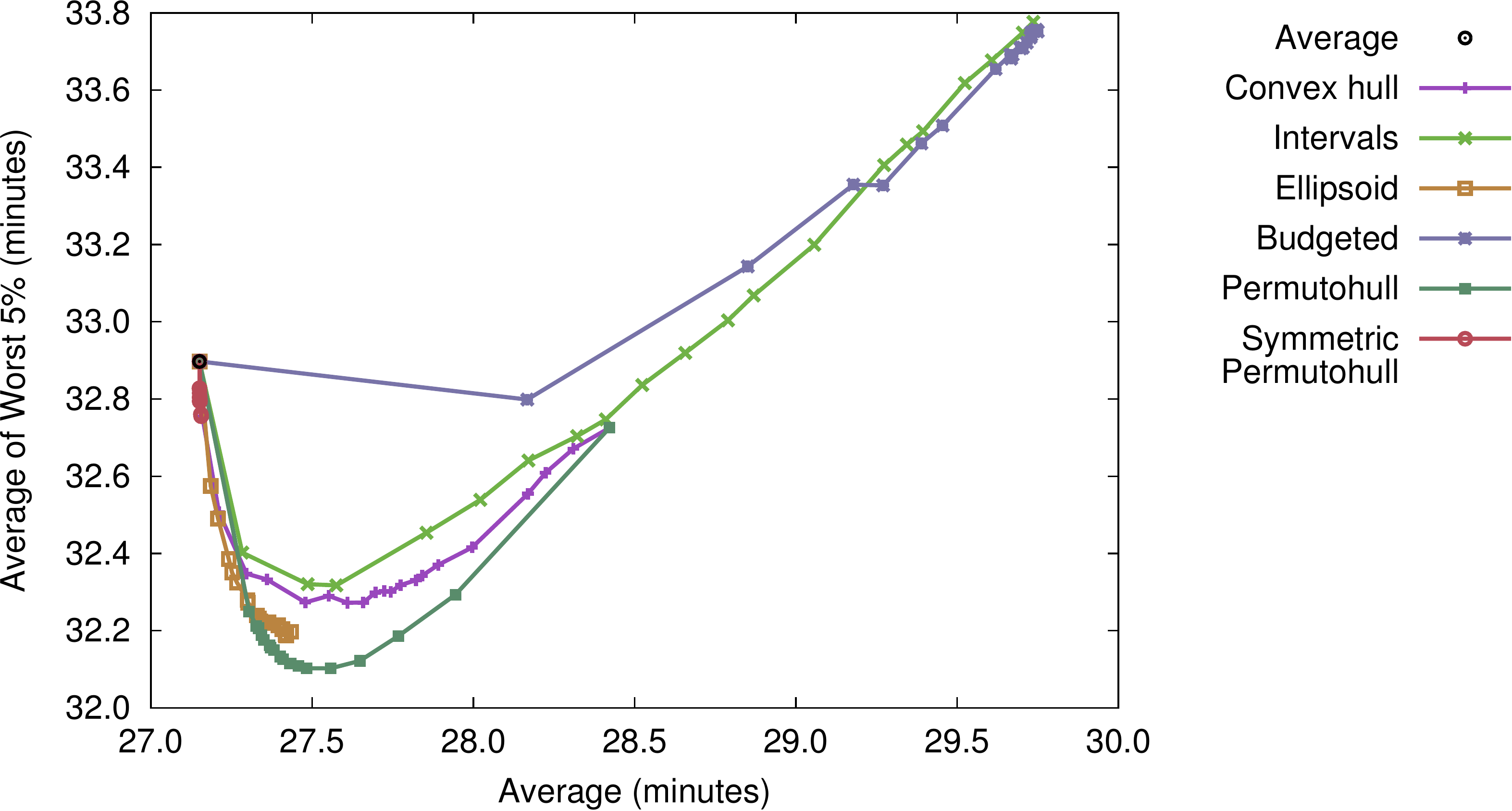}}      
    \hfill
\subfigure[out-sample]{\label{av-cvar-out-evenings}\includegraphics[width=\textwidth]{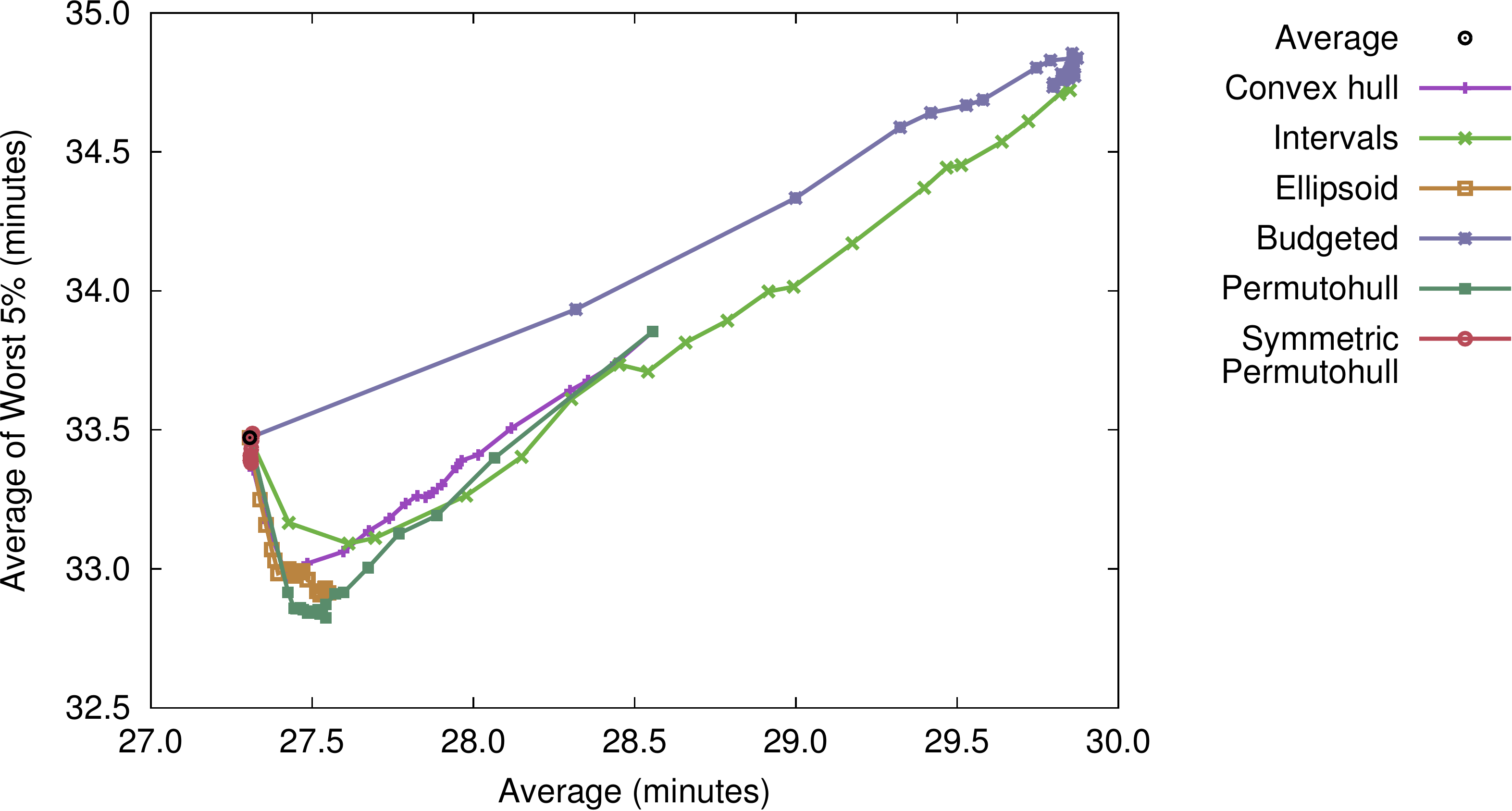}}      
\caption{Average vs CVaR performance, evenings dataset.}\label{av-cvar-evenings}
\end{figure}

% % % % % % % % % % % % 

\begin{figure}[htbp]
\centering
\subfigure[in-sample]{\label{av-max-in-tuesdays}\includegraphics[width=\textwidth]{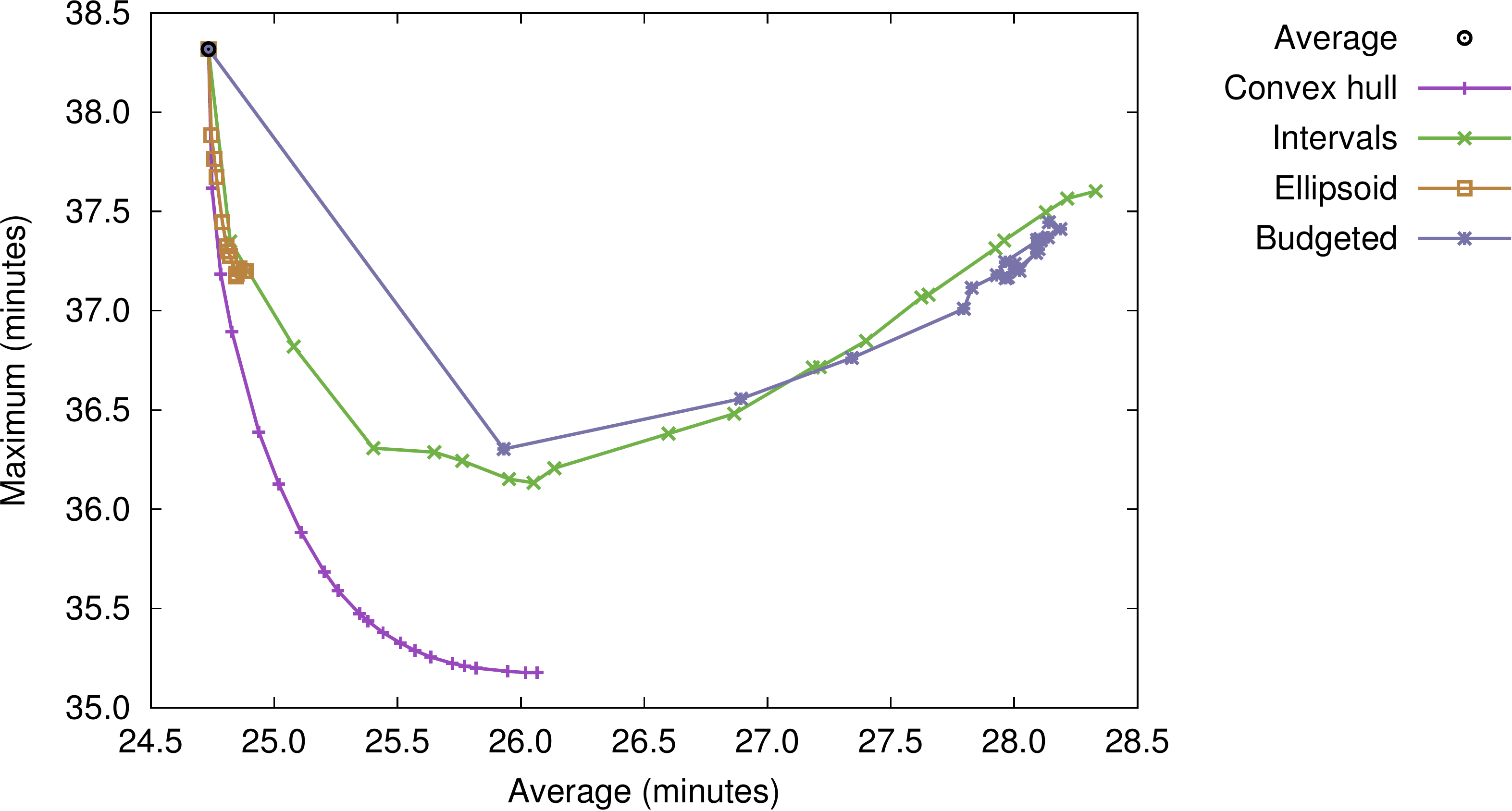}}  
    \hfill
\subfigure[out-sample]{\label{av-max-out-tuesdays}\includegraphics[width=\textwidth]{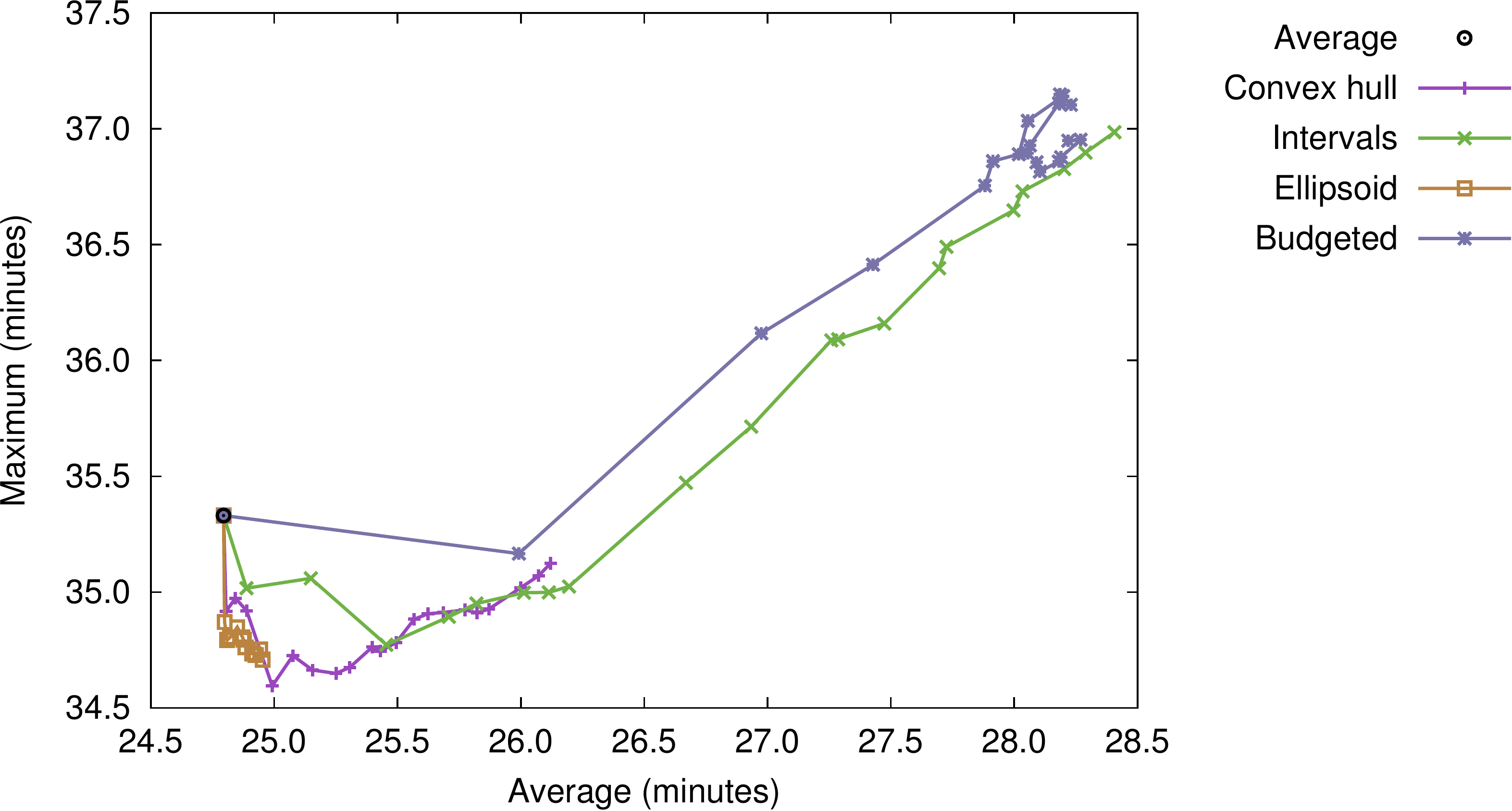}}      
\caption{Average vs worst-case performance, Tuesdays dataset.}\label{av-max-tuesdays}
\end{figure}

\begin{figure}[htbp]
\centering
\subfigure[in-sample]{\label{av-cvar-in-tuesdays}\includegraphics[width=\textwidth]{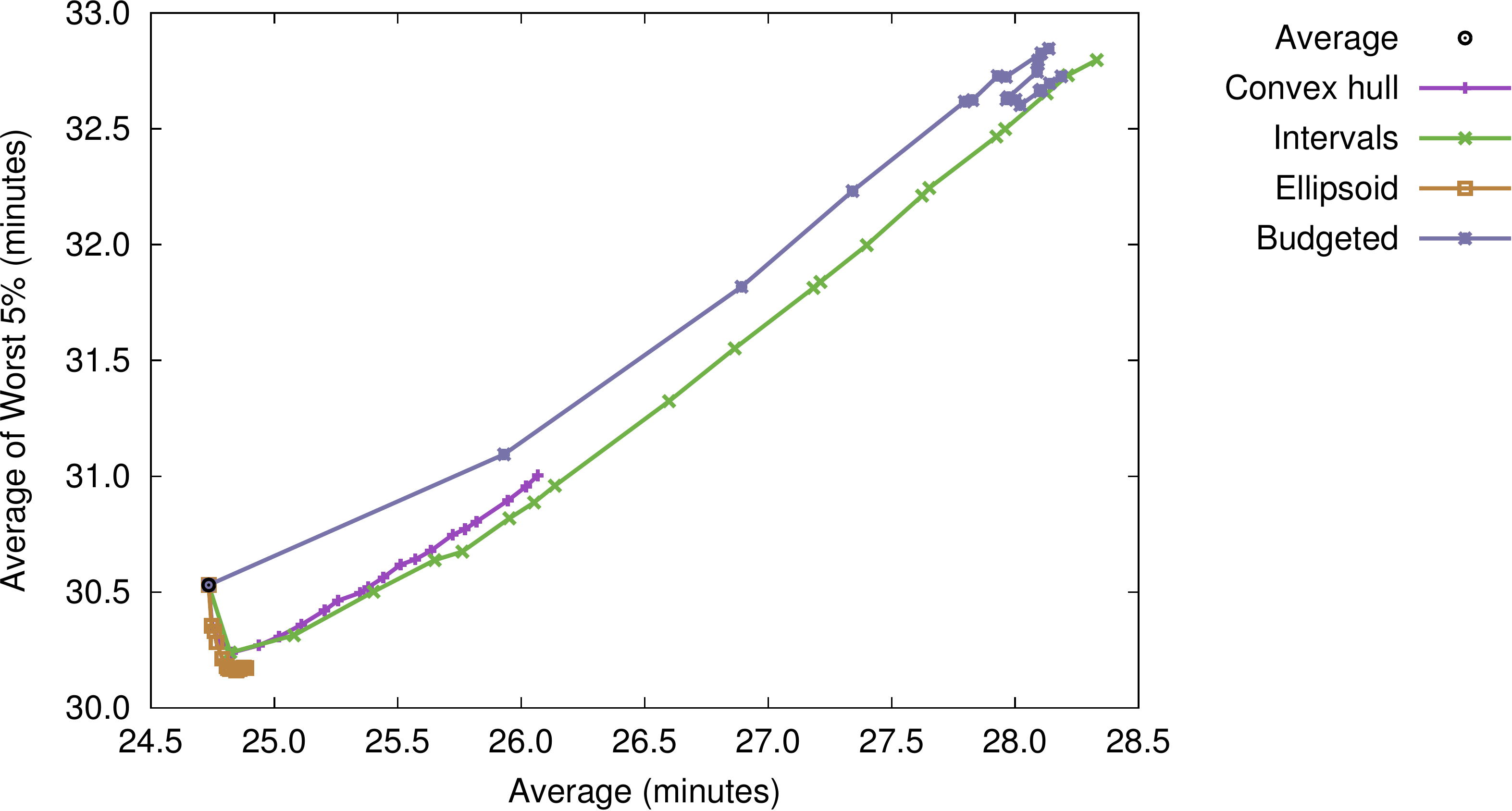}}      
    \hfill
\subfigure[out-sample]{\label{av-cvar-out-tuesdays}\includegraphics[width=\textwidth]{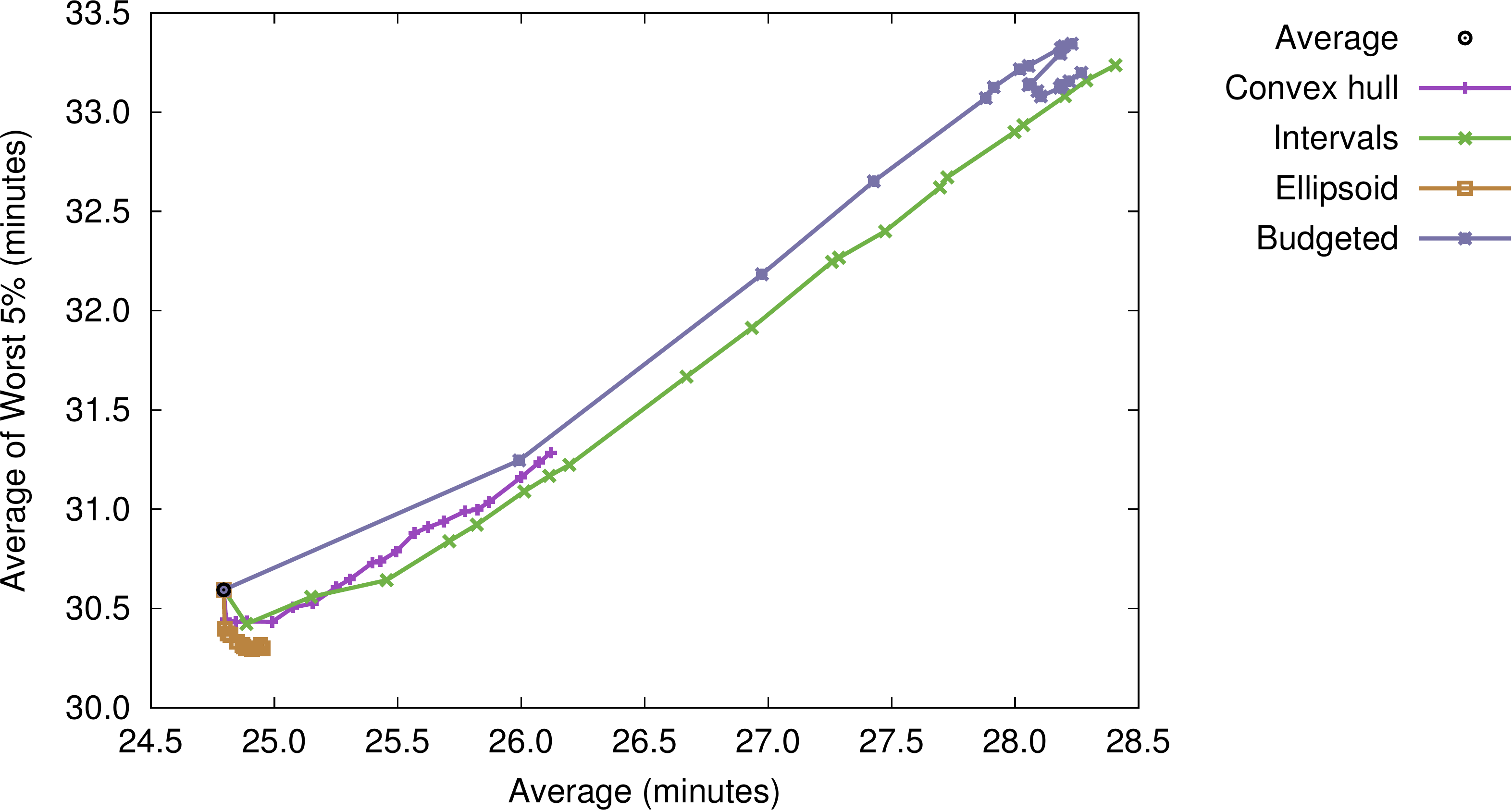}}      
\caption{Average vs CVaR performance, Tuesdays dataset.}\label{av-cvar-tuesdays}
\end{figure}

% % % % % % % % % % % % 

\begin{figure}[htbp]
\centering
\subfigure[in-sample]{\label{av-max-in-weekends}\includegraphics[width=\textwidth]{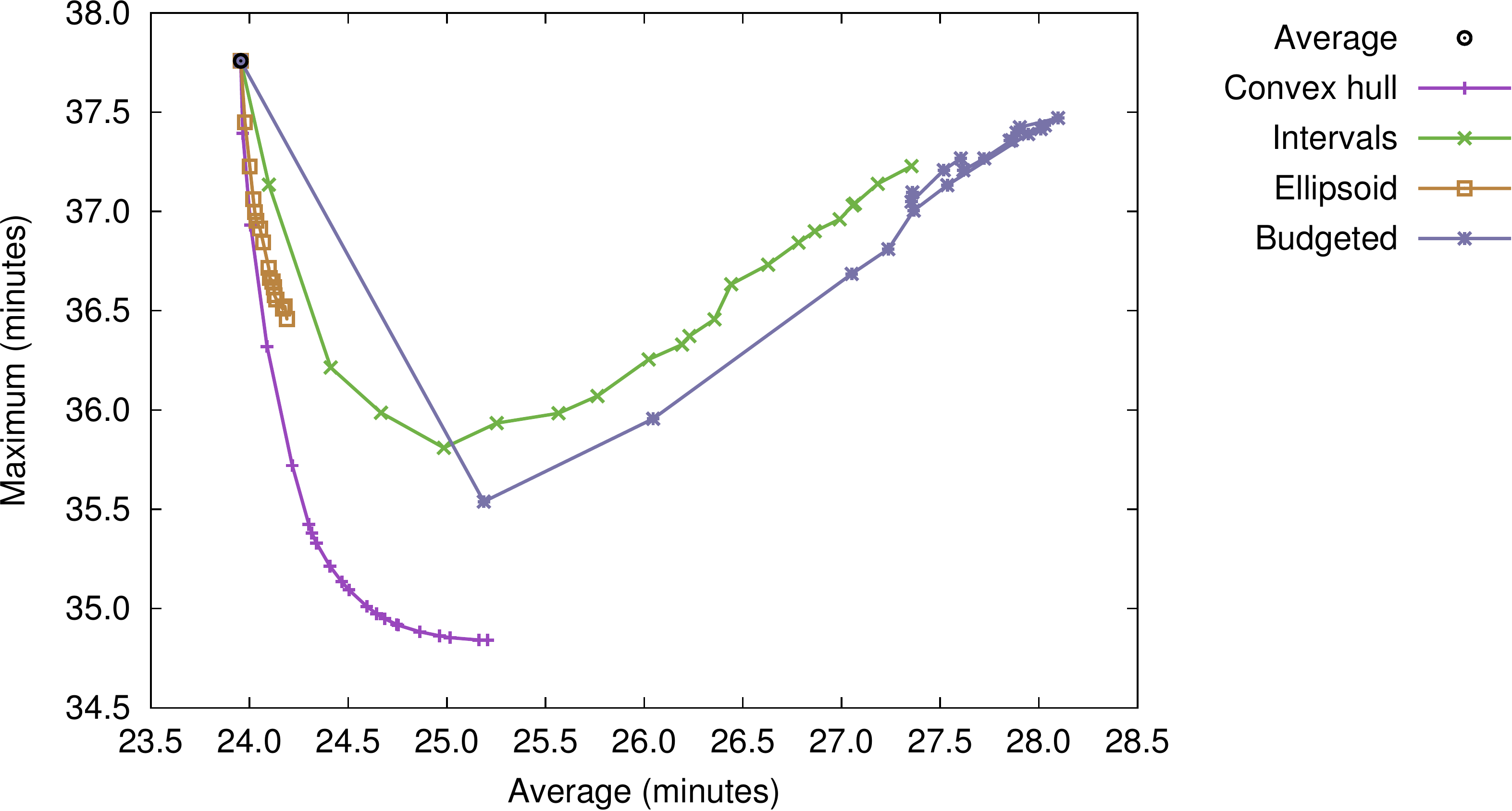}}  
    \hfill
\subfigure[out-sample]{\label{av-max-out-weekends}\includegraphics[width=\textwidth]{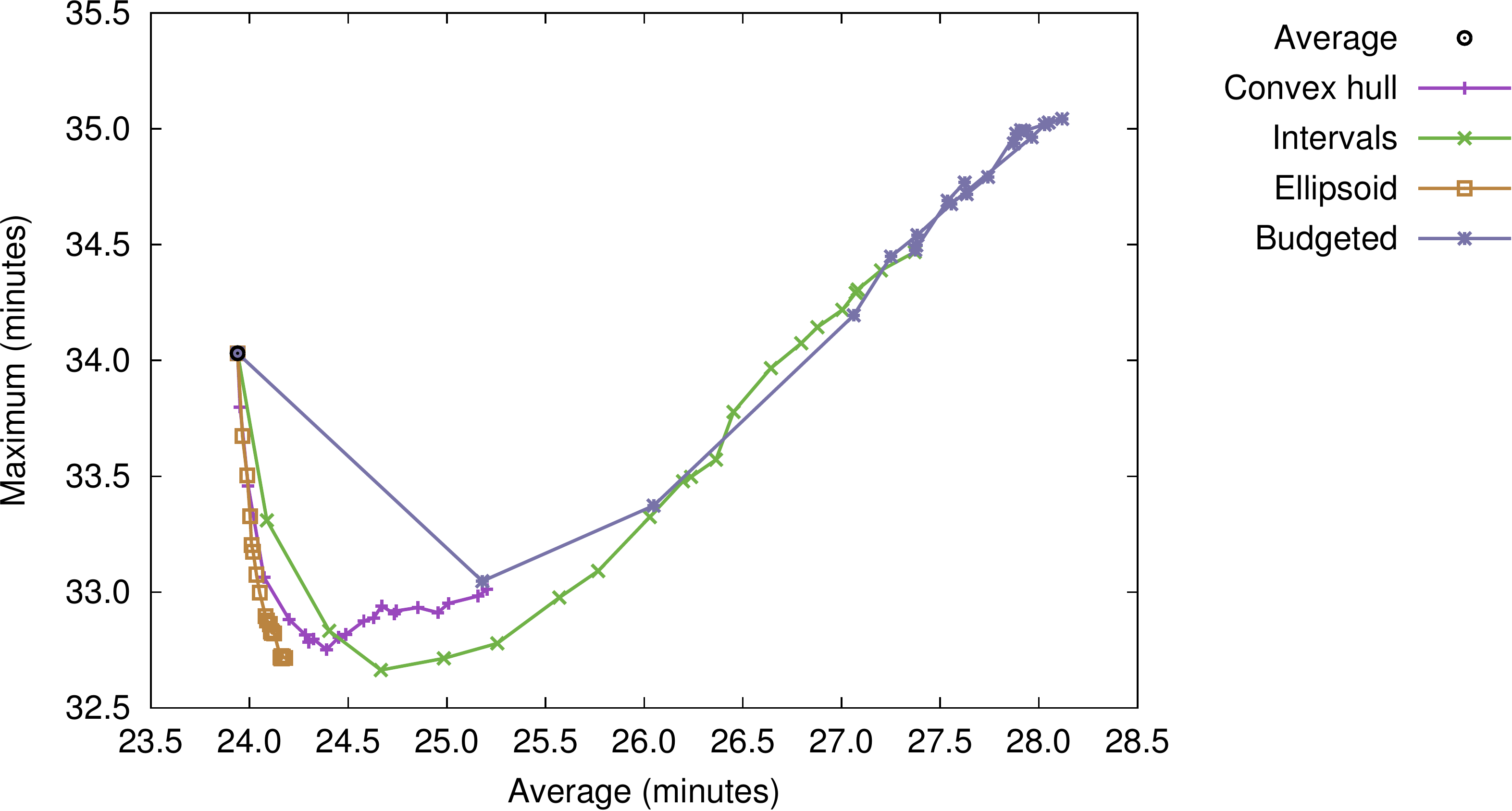}}      
\caption{Average vs worst-case performance, weekends dataset.}\label{av-max-weekends}
\end{figure}

\begin{figure}[htbp]
\centering
\subfigure[in-sample]{\label{av-cvar-in-weekends}\includegraphics[width=\textwidth]{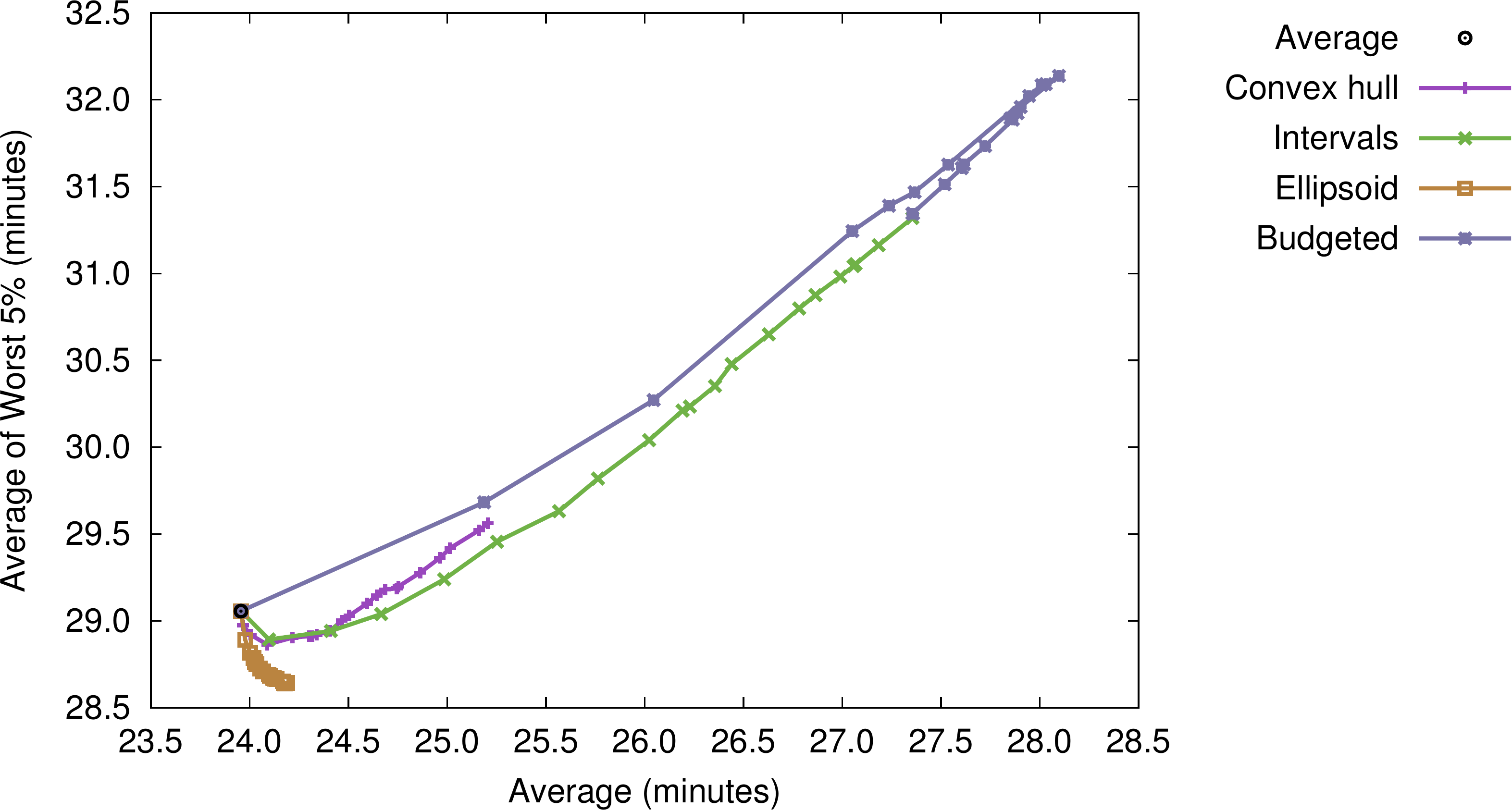}}      
    \hfill
\subfigure[out-sample]{\label{av-cvar-out-weekends}\includegraphics[width=\textwidth]{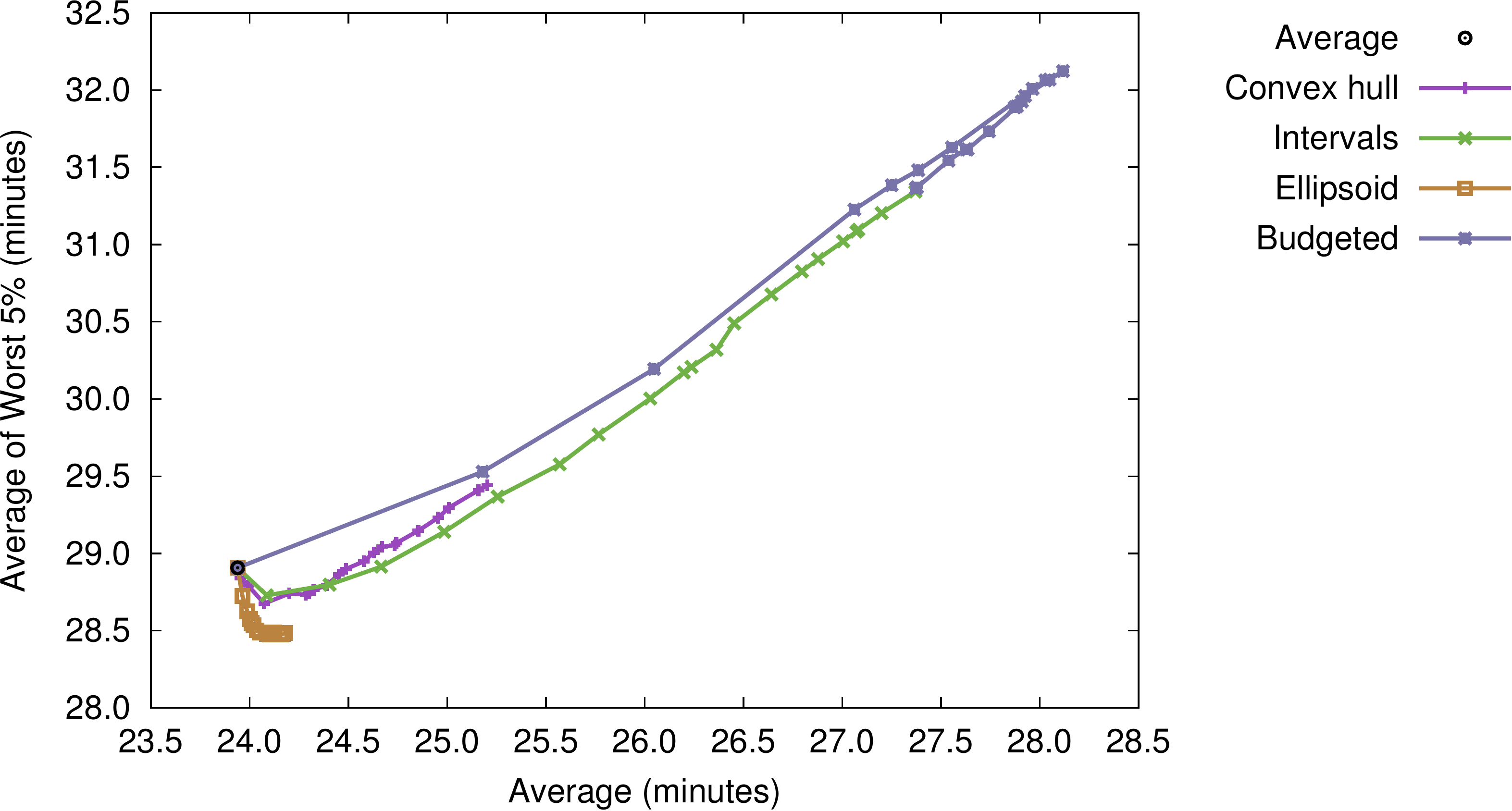}}      
\caption{Average vs CVaR performance, weekends dataset.}\label{av-cvar-weekends}
\end{figure}

% % % % % % % % % % % % 

\begin{figure}[htbp]
\centering
\subfigure[in-sample]{\label{av-max-in-all}\includegraphics[width=\textwidth]{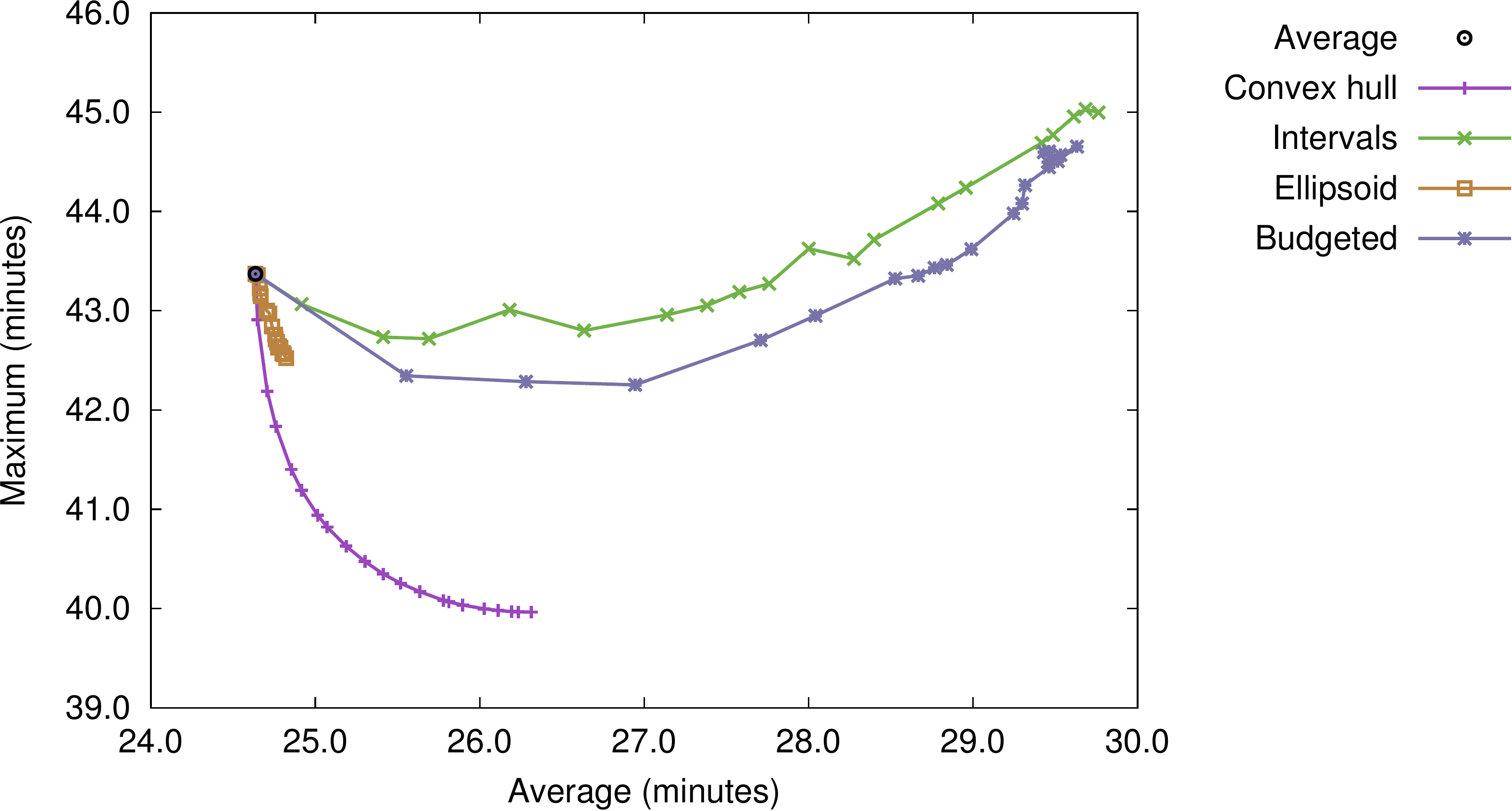}}  
    \hfill
\subfigure[out-sample]{\label{av-max-out-all}\includegraphics[width=\textwidth]{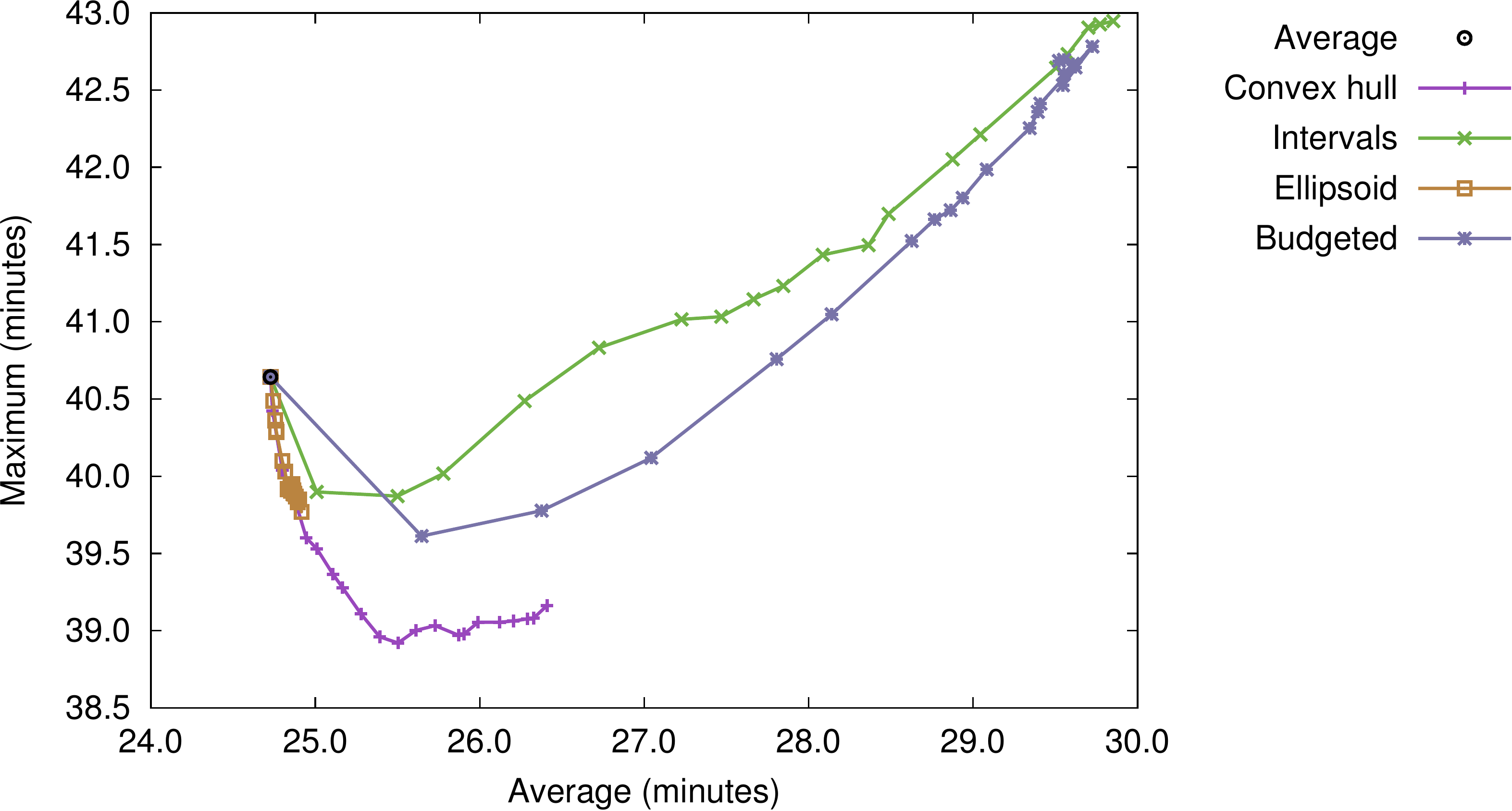}}      
\caption{Average vs worst-case performance, complete dataset.}\label{av-max-all}
\end{figure}

\begin{figure}[htbp]
\centering
\subfigure[in-sample]{\label{av-cvar-in-all}\includegraphics[width=\textwidth]{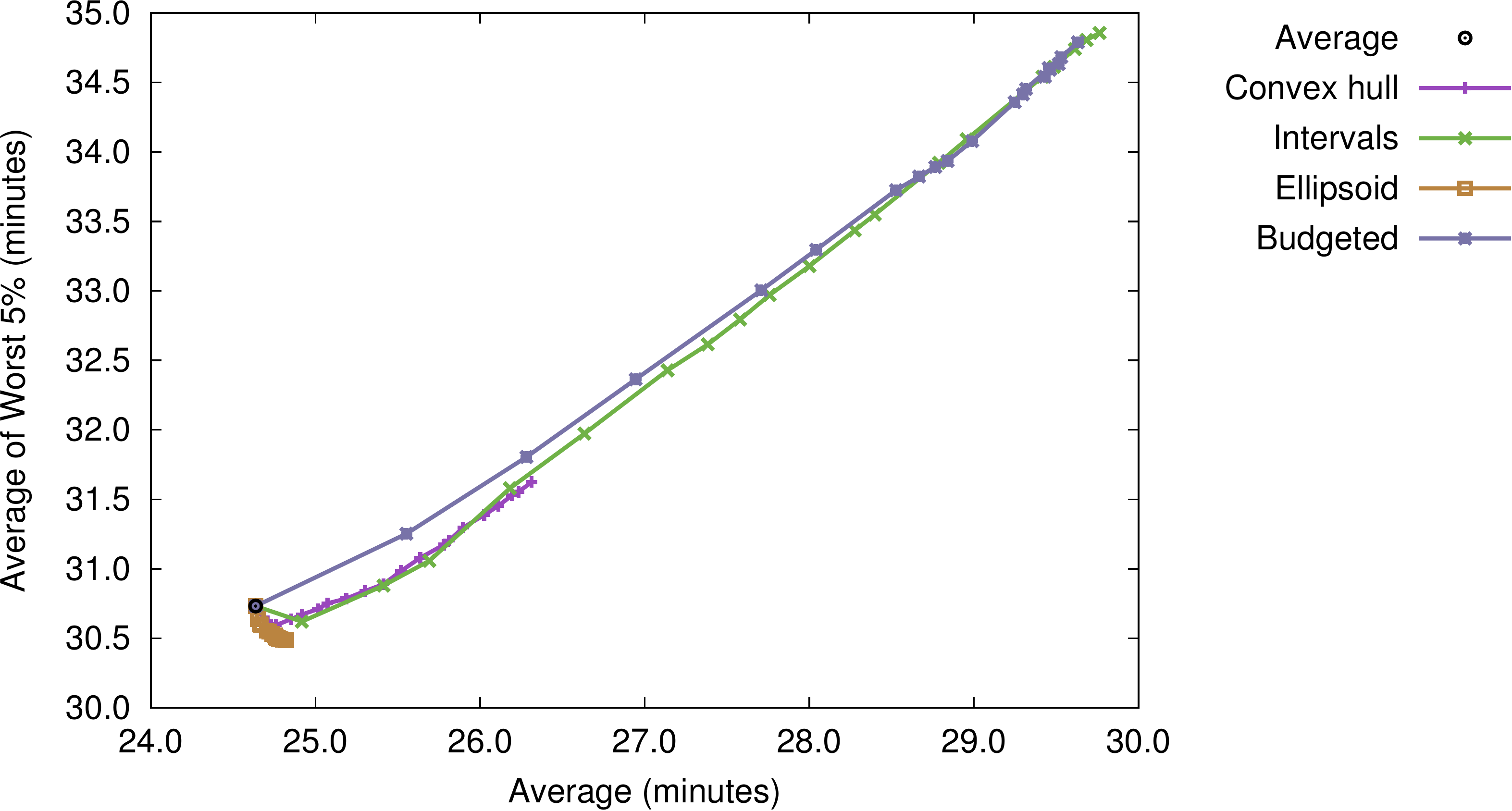}}      
    \hfill
\subfigure[out-sample]{\label{av-cvar-out-all}\includegraphics[width=\textwidth]{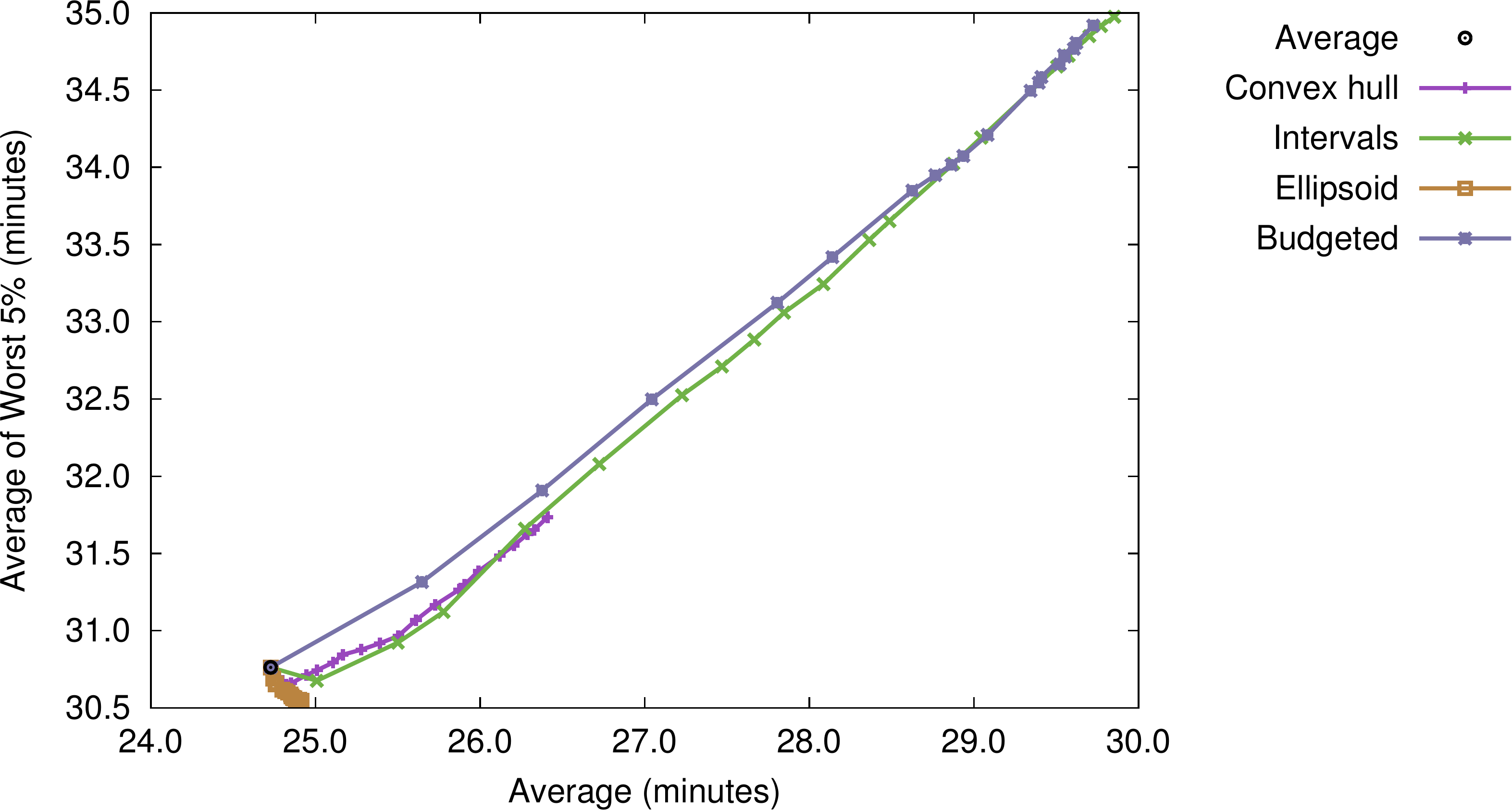}}      
\caption{Average vs CVaR performance, complete dataset.}\label{av-cvar-all}
\end{figure}

\end{document}